\documentclass{article}
\usepackage{isolatin1}
\usepackage{amsmath}
\usepackage{amsthm}
\usepackage{amssymb}

\usepackage[all,dvips]{xy}
\usepackage{diagrams}
\usepackage{comment}

\excludecomment{longversion}
\includecomment{shortversion}

\usepackage{ifthen}

\theoremstyle{plain}
\newtheorem{Theorem}{Theorem}[section]
\newtheorem{Proposition}[Theorem]{Proposition}
\newtheorem{Corollary}[Theorem]{Corollary}
\newtheorem{Lemma}[Theorem]{Lemma}

\theoremstyle{definition}
\newtheorem{Definition}[Theorem]{Definition}

\newtheorem{Fact}[Theorem]{Fact}

\theoremstyle{remark}

\newcommand{\under}{\mbox{$\setminus\!\!\setminus$} }
\newcommand{\inversions}{\mathcal{I}}
\newcommand{\set}[1]{\{\,#1\,\}}
\newcommand{\setof}[1]{[#1]}
\newcommand{\hspp}{\hspace{5mm}}
\newcommand{\vspp}{\vspace{12pt}}

\newdir{ >}{{}*!/-5pt/\dir{>}}
\newcommand{\LP}{\mathcal{L}}
\newcommand{\Paths}{\LP}
\newcommand{\Perm}{\mathcal{P}erm}
\newcommand{\dreftrclosure}{\unlhd}
\newcommand{\dep}{\,{\unlhd}^{\bullet}\,}
\newcommand{\dtrclosure}{\lhd}
\newcommand{\jmup}{\nearrow}
\newcommand{\mjdown}{\searrow}

\newcommand{\N}{\mathbb{N}}
\newcommand{\covers}{\prec}
\newlength{\jimilen}

\newlength{\jimilend}

\newcommand{\ji}[1]{%
  {}_{\vee}%
  \settowidth{\jimilend}{$#1$}%
  \ifdim\jimilend<\jimilen%
  {{#1}}%
  \else%
  {(#1)}%
  \fi%
}
\newcommand{\mi}[1]{%
  {}^{\wedge}
  \settowidth{\jimilend}{$#1$}%
  \ifdim\jimilend<\jimilen%
  {{#1}}%
  \else%
  {(#1)}%
  \fi%
}

\renewcommand{\vec}[1]{x_{#1}}
\newcommand{\minji}{\textstyle\min_{\vee}}
\newcommand{\maxji}{\textstyle\max_{\vee}}
\newcommand{\minmi}{\textstyle\min_{\land}}
\newcommand{\maxmi}{\textstyle\max_{\land}}

\newgraphescape{E}[2]{[]*+{#2}="#1"*+{{\mbox{\hspace{12mm}}}^{#1}}"#1"}

\newgraphescape{R}[2]{
  []*+{#1}="1"-[d]*+{{#1}_{\star}}="2"
  -[ld]*+{#2}="3"-[d]*+{{#2}_{\star}}="4"
  -[r(2)u]*+{\kappa(#2)}="5"
  -[u(3)]*+{\kappa(#2)^{\star}}="6"-"1"
}

\newgraphescape{L}[2]{
  []*+{#1}="1"-[d]*+{{#1}_{\star}}="2"-[rd]*+{#2}="3"
  -[d]*+{{#2}_{\star}}="4"
  -[l(2)u]*+{\kappa(#2)}="5"
  -[u(3)]*+{\kappa(#2)^{\star}}="6"-"1"
}

\newgraphescape{M}{
  []="0"(
  :@{.}[r(3)],
  [u](:@{.}[r(3)])[u](:@{.}[r(3)])[u]:@{.}[r(3)],
  :@{.}[u(3)]="N",
  [r](:@{.}[u(3)])[r](:@{.}[u(3)])[r](:@{.}[u(3)]="NE")
  )
}

\newgraphescape{A}{
  []:@{.}[r(2)]
}
\newgraphescape{B}{
  []:@{.}[r(1.5)u(0.7)]
}
\newgraphescape{C}{
  []:@{.}[u(2)]
}
\newgraphescape{N}{
  []="0"(
  !A!B!C,
  !A!C="acb"!B,
  !B="bac"!A!C,
  !B!C="bca"!A,
  !C="cab"!A!B,
  !C!B!A
  )
}

\newgraphescape{D}[2]{
  [d(#1)r(#2)]*+{\bullet}
}

\newgraphescape{Q}[6]{
  []*+{#4}="S"="0"
  -[ru]*+{#1}="1"
  -[u]*+{#2}="2"
  -[lu]*+{#5}="4"
  -[d(1.5)l(1)]*+{#3}="3"
  -"S"
  [u(1.5)]*+{#6}
  "S"
}
\newgraphescape{P}[6]{
  []*+{#4}="S"="0"
  -[lu]*+{#1}="1"
  -[u]*+{#2}="2"
  -[ru]*+{#5}="4"
  -[d(1.5)r(1)]*+{#3}="3"
  -"S"
  [u(1.5)]*+{#6}
  "S"
}

\newcommand{\mygraph}[2][3em]{%
  \mbox{%
    $\begin{array}{@{\hspace{0mm}}c@{\hspace{0mm}}}
      \xy\xygraph{!{<#1,0cm>:<0cm,#1>::}#2}\endxy%
    \end{array}$%
  }%
}

\author{Luigi Santocanale\\%
  \texttt{lsantoca@cmi.univ-mrs.fr}\\%
  LIF/CMI, Université de Marseille}
\title{Congruences of Multinomial Lattices}

\begin{document}

\maketitle

\begin{abstract}
  We study the congruence lattices of the multinomial lattices
  $\LP(v)$ introduced by Bennett and Birkhoff \cite{bb}. Our main
  motivation is to investigate Parikh equivalence relations that model
  concurrent computation. We accomplish this goal by providing an
  explicit description of the join dependency relation between two
  join irreducible elements and of its reflexive transitive closure.
  The explicit description emphasizes several properties and makes it
  possible to separate the equational theories of multinomial lattices
  by their dimensions.
  
  In their covering of non modular varieties \cite{jipsen} Jipsen
  and Rose define a sequence of equations \ref{eq:sdn}, for $n \geq
  0$. Our main result sounds as follows: \emph{if $v = (v_{1},\ldots
    ,v_{n}) \in \N^{n}$ and $v_{i} > 0$ for $i = 1,\ldots ,n$, then
    the multinomial lattice $\LP(v)$ satisfies $SD_{n-1}(\land)$ and
    fails $SD_{n-2}(\land)$}.
\end{abstract}

%%% Local Variables: 
%%% mode: latex
%%% TeX-master: "0"
%%% End: 
 
\section*{Introduction}

Multinomial lattices were introduced in \cite{bb} in the context of an
order theoretic investigation of rewrite systems associated to common
algebraic laws. From this point of view, they form an order theoretic
counterpart of the commutativity law.

As a family of finite lattices, multinomial lattices subsume two well
known families. The binomial lattices $\LP(p)$, $p \in \N^{2}$, are
also known as lattices of lattice-paths, since their elements are
paths in the discrete plane from $0$ to $p$.  Counting properties of
paths in the set $\LP(p)$ have been intensively investigated, see
\cite{mohanty,krattenthaler}.  Order theoretic properties of $\LP(p)$
have been studied in \cite{bb,pinzani}.  Among these properties, these
lattices are distributive.

The second family of lattices are the permutoedra $\Perm(n)$, $n \geq
0$. Elements of $\Perm(n)$ are permutations on the set $\set{1,\ldots
  ,n}$.  It was shown in \cite{guilbaud} that this set, endowed with
the weak Bruhat order, is a lattice, a result later generalized to all
finite Coxeter groups \cite{bjorner}.  The lattice structure of
$\Perm(n)$ has been deeply investigated as well, see
\cite{yanagimoto,duquenne,markowsky,poly-barbut1,poly-barbut2,caspard}.

A multinomial lattice $\LP(v)$ has as underlying set the collection of
all ``discrete'' paths from $0$ to $v$, where $v = (v_{1},\ldots
,v_{n})$ is a vector in $\N^{n}$ and $n$, the dimension, can be an
arbitrary positive integer.  The paths we consider are discrete in
that they add $1$ to just one coordinate at each time unit.  For this
reason we used to refer to multinomial lattices as lattices of paths in
higher dimension.  These paths are in bijection with words $w$ over an
alphabet $\Sigma = \set{a_{1},\ldots ,a_{n}}$ such that the number of
letters $a_{i}$ occurring in $w$ is equal to $v_{i}$.

It was shown in \cite{bb} that $\LP(v)$, as se set, can be endowed
with an order structure which turns out to be a lattice structure. If
$n$, the dimension, is $2$, then $\LP(v)$ is the usual lattice of
lattice paths. Let $1^{n}$ be the vector in $\N^{n}$ with just $1$'s
at each coordinate, then $\LP(1^{n})$ is (order isomorphic to) the
permutoedron $\Perm(n)$.  In this sense, the multinomial lattices are
a common generalization of the permutoedra and of lattices of lattice
paths.

A main motivation for us to approach multinomial lattices has been
investigating Parikh equivalence relations that arise when modeling
concurrent computation. If $\Sigma$ is an alphabet, an equivalence
relation $\sim$ on $\Sigma^{n}$ is Parikh if $w \sim u$ implies that
the number of occurrences of a letter $\sigma$ in $w$ and in $u$ are
the same, for each $\sigma \in \Sigma$.  Lattice congruences of
$\LP(v)$ give rise to a class of Parikh's equivalence relations which
enjoy a property of interest in concurrency: if $w \sim u$, then we
can find a sequence $w=x_{0},x_{1},\ldots ,x_{n} = u$ such that $x_{i
  + 1}$ is obtained from $x_{i}$ by switching two contiguous letters.
Consequently, our goal has been to understand the congruence lattice
of the multinomial lattices $\LP(v)$.

We accomplish this goal by providing an \emph{explicit description of
  the join dependency relation} among join irreducible elements. The
explicit description emphasizes several properties of the lattices
$\LP(v)$'s, for example these lattices are bounded. Among the
properties a curious one: a sequence of join irreducible elements
related by the join dependency relation cannot have length greater
than $n-2$, if $v \in \N^{n}$. This property suffices to separate the
equational theories of the $\LP(v)$'s by dimension.  In \cite[\S
4.2]{jipsen} a family of equations \ref{eq:sdn}, $n \geq 0$, is
introduced, that can be taken as a measure of meet semidistributivity
of
a finite lattice.\footnote{%
  Analogous considerations and results holds for join
  semidistributivity and the dual equations $SD_{n}(\vee)$.  } Our
main result can be phrased as follows: \emph{if $v = (v_{1},\ldots
  ,v_{n}) \in \N^{n}$ and $v_{i} > 0$ for $i = 1,\ldots ,n$, then the
  multinomial lattice $\LP(v)$ satisfies $SD_{n-1}(\land)$ and fails
  $SD_{n-2}(\land)$.}  For example, assuming that dimensions are not
degenerate, lattices of lattice-paths are distributive and not reduced
to a point, a lattice $\LP(v)$ with $v \in \N^{3}$ is neardistributive
but not distributive, and so on.

While the results presented here introduce a notion of dimension for
multinomial lattices, it doesn't appear to exist %at the present moment
a direct relation with the order dimension, determined for multinomial
lattices in \cite{flath}.

\newpage
\tableofcontents

%%% Local Variables: 
%%% mode: latex
%%% TeX-master: "0"
%%% ispell-dictionary: "english"
%%% End: 

\section{Notation and Background}

\subsection{Words and functions}

With $\setof{n}$ we shall denote the set $\set{1,\ldots ,n}$. Recall
that a word over an alphabet $\Sigma$ is a function $w : \setof{n}
\rTo \Sigma$, $n$ being the length of the word $w$. The notations
$w(i)$ and $w_{i}$ are therefore equivalent. Similarly, for a
permutation $\sigma : \setof{n} \rTo \setof{n}$, $\sigma(i) =
\sigma_{i}$.  With $(i,j)$ we shall denote the permutation sending $i$
to $j$ and viceversa, and fixing all the other members of $\setof{n}$,
and  $\sigma^{i}$ will denote the permutation $(i, i+ 1)$.

\subsection{Congruences of  finite lattices}

A lattice is an ordered set with the property that every finite
non-empty subset has a least upper bound and a greatest lower bound,
see the standard literature on lattices \cite{Birk,BM,gratzer,davey}.
Let $x \vee y$ and $x\wedge y$ denote the least upper bound and the
greatest lower bound of the finite set $\set{x,y}$, respectively. With
respect to the binary operations $\vee$ and $\wedge$, lattices are
algebraic structures, also meaning that the order is equationally
definable and determined by the two binary operations.

In this paper we shall be studying finite lattices and, when
considering a lattice, we shall assume it is finite unless otherwise
stated.  It is a standard argument that a finite lattice has a top and
a bottom elements. All the structure of a lattice is determined by the
order relation restricted to join irreducible elements and meet
irreducible elements, see \cite{wille}.  An element $j$ is join
irreducible if $j = x \vee y$ implies $j = x$ or $j = y$, and an
analogous property defines a meet irreducible element.  There is an
order theoretic characterization of being join/meet irreducible. To
this goal, recall that -- for an arbitrary ordered set -- $x < y$ is a
cover (noted $x \covers y$) if the closed interval $[x,y]$ has only
two elements; then we say that $x$ is a lower cover of $y$ and $y$ is
an upper cover of $x$. An element $j$ of a lattice is join irreducible
if and only if it has a unique lower cover, which is denoted by
$j_{\star}$; $m$ is meet irreducible if it has a unique upper cover
$m^{\star}$.

We recall the definitions/characterizations of standard relations
between join/meet, meet/join  irreducible elements:
%% \begin{itemize}
%% \item $j \jmup m$ iff $j \not\leq m$ and $j \leq m^{*}$,
%% \item $m \mjdown j$ iff $j \not\leq m$ and $j_{*} \leq m$,
%% \item $j D j'$ iff $j \neq j'$ and $j \jmup m \mjdown j$ for
%%   some meet irreducible $m$.
%% \end{itemize}
\begin{align}
  \label{eq:up}
  j \jmup m & \text{ iff } j \not\leq m \text{ and } j \leq m^{*}\,,
  \\
  \label{eq:down}
  m \mjdown j & \text{ iff } j \not\leq m \text{ and } j_{\star}
  \leq m\,,\\
  \intertext{and of the join dependency relation $D$ between join
    irreducible elements:}
  \label{eq:D}
  j D j' & \text{ iff } j \neq j', j \jmup m \mjdown j' \text{ for some meet
  irreducible } m\,.
\end{align}
The meet dependency relation is defined as expected: $mD^{d}m'$ iff $m
\neq m'$, $m \mjdown j \jmup m'$ for some $j$.  These relations are
central in the theory of finite lattices, as we explain next. Since a
lattice is an algebraic structure, we can define a congruence on a
lattice $L$ as an equivalence relation $\theta \subseteq L \times L$
compatible with the lattice operations, i.e.  such that $x \theta y$
implies $x \vee z \theta y \vee z$ and $x \land z \theta y \land z$.
The quotient $L/\theta$ is then a lattice and the canonical projection
is a lattice homomorphism. The following Proposition has motivated us
to study the join dependency relation $D$ and its transitive closure
(reflexive and transitive closure) $\dtrclosure$ (resp.
$\dreftrclosure$) in multinomial lattices.
\begin{Proposition}
  The congruences of a lattice $L$ are in bijection with subsets $S$
  of join irreducible elements that are closed under the $D$-relation.
  The bijection is an order reversing isomorphism of lattices.
\end{Proposition}
For a proof, the reader may consult \cite[\S 2.34 and \S
11.10]{freese}. It is convenient to explicit the bijection. Given a
congruence $\theta$, the $S_{\theta}$ is defined as the set of join
irreducible elements $j$ such that $j\theta j^{\star}$ does not hold,
i.e. that are not collapsed with their lower cover under the
congruence $\theta$. The latter can be recovered from $\theta$, since

\begin{align}
  \label{eq:congrfromS}
  x \theta y \text{ if and only if } &
  \set{z \in S_{\theta}\,|\,z\leq x}
  = \set{z \in S_{\theta}\,|\,z\leq y}\,.
\end{align}
A lattice $L$ is \emph{semidistributive} if the
conditions
  \begin{align}
  \label{eq:sd}
  \tag{$SD(\land)$}
  x \land y  = x \land z & \Rightarrow x \land (y \vee z) = x\land y
  \\
  \label{eq:sdvee}
  \tag{$SD(\vee)$}
  x \vee y  = x \vee z & \Rightarrow x \vee (y \land z) = x\vee y
\end{align}
hold in $L$. Equivalently, a lattice is semidistributive if for each
join irreducible there exists a unique $m$ such that $j\jmup m \mjdown
j$ and the dual condition hold, see \cite[\S 2.6]{freese}.  We shall
write $\kappa(j)$ for such an $m$ (and $\kappa^{d}(m)$, dually).

A lattice is \emph{bounded} if it is the quotient of a finitely
generated free lattice (with top and bottom which usually is
infinite), and the quotient map has both a left and a right adjoints.
There are several characterizations of the notion of bounded lattice,
see in particular \cite{day}. Among them we shall use the following
one, see \cite[\S 2.53]{freese}:
\begin{Lemma}
  \label{lemma:bounded}
  A lattice is bounded if and only if it is semidistributive and
  the join dependency relation $D$ contains no cycle.
\end{Lemma}

\subsection{Lattices of lattice paths, i.e. binomial lattices}
\label{sec:latticepaths}
We shall denote by $\LP(n,m)$ the set of words $w$ over the alphabet
$\Sigma = \set{a,b}$ such that $|w|_{a} = n$ and $|w|_{b} = m$. We
represent these words as paths in the 2-dimensional space from $(0,0)$
to $(n,m)$: if $w \in \LP(n,M)$, then the path $f_{w}: \set{0,\ldots ,
  n +m} \rTo \N \times \N$ is defined by induction as follows: $f(0) =
(0,0)$, $f(i) = f(i-1) + (1,0)$ if $w_{i} = a$ and $f(i) = f(i-1) +
(0,1)$ if $w_{i} = b$, for $i = 1,\ldots ,n+m$.  The following diagram
represents the word/path $abaabb \in \LP(3,3)$:
$$
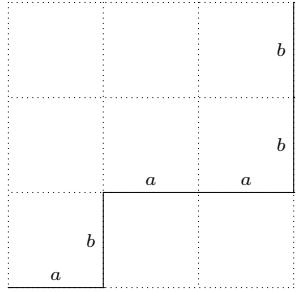
\xygraph{
  []!M
  "0"-[r]^{a}-[u]^{b}
  -[r]^{a}-[r]^{a}
  -[u]^{b}-[u]^{b}
}
$$
The rewrite relation $\rightarrow$ on $\LP(n,m)$ is defined by:
\begin{align*}
  w \rightarrow u &
  \text{ iff } w = w_{1}abw_{2}
  \text{ and }
  u = w_{1}baw_{2}\,.
\end{align*}
This terminating and confluent rewrite system gives rise to an order
relation $\leq$ which is a distributive lattice. Indeed, we have $f
\leq g$ if and only if $f(i) \leq g(i)$ for $i \in \set{0,\ldots ,n +
  m}$, where the pointwise order is defined as follows: $(x_{1},y_{1})
\leq (x_{2},y_{2})$ if and only if $x_{2} \leq x_{1}$ and $y_{1} \leq
y_{2}$.

A join irreducible of $\LP(n,m)$ is a word of the form
$a^{x}b^{y}a^{\bar{x}}b^{\bar{y}}$, where $(\bar{x},\bar{y}) = (n,m) -
(x,y)$, $0 <\bar{x} \leq n$, and $0 < y \leq m$. Clearly the latter
two conditions are equivalent to $0\leq x < n$ and $0 < y \leq m$. A
join irreducible is therefore uniquely determined by a vector $(x,y)$
with this property, and we shall use the notation for $\ji{x,y}$ for
the join irreducible $a^{x}b^{y}a^{\bar{x}}b^{\bar{y}}$. In
\cite{krattenthaler} a join irreducible is called a turn to NorthEast
(a NE turn). The fact that a path can be described uniquely by its
NE-turns corresponds to the lattice theoretic property that an element
of a finite lattice is the join of the join irreducible elements below
it, among which we can retain the antichain of maximal elements.

Similarly, a meet irreducible has the form
$b^{y}a^{x}b^{\bar{y}}a^{\bar{x}}$ with $0 < x \leq n$ and $0 \leq y <
m$ and  is uniquely determined by the vector $(x,y)$: we shall use
the notation $\mi{x,y}$ for $b^{y}a^{x}b^{\bar{y}}a^{\bar{x}}$.
Observe that
\begin{align}
  \label{eq:joinirrpaths}
  \ji{x,y} \leq \mi{z,w}
  & \text{ iff } z \leq x \text{ or } y \leq w\,.
\end{align}
Since $\LP(n,m)$ is distributive, it is semidistributive as well.
Recalling the definitions \eqref{eq:up} and \eqref{eq:down} of the
relations $\jmup$ and $\mjdown$, we observe that in $\LP(n,m)$
\begin{align}
  \notag
  \ji{x,y} \jmup \mi{z,w}
  & 
  \text{ iff } \mi{z,w} \mjdown \ji{x,y}%
  \\
  \label{eq:jmupllp}
  & \text{ iff } z = x + 1 \text{ and } w = y - 1\,.
\end{align}
For example, $\kappa(aabbab)= (baaabb)$ or 
$\kappa \ji{2,2} = \mi{3,1}$. The join irreducible $\ji{2,2}$ and
$\kappa \ji{2,2}$ are represented in following diagram:
$$
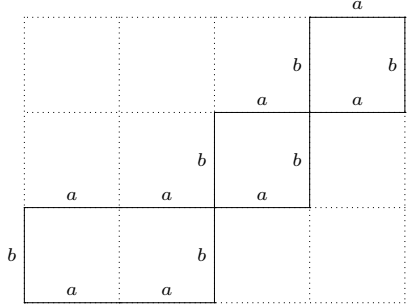
\xygraph{
  []="0"(
  :@{.}[r(4)],
  [u](:@{.}[r(4)])
  [u](:@{.}[r(4)])
  [u](:@{.}[r(4)]),
  :@{.}[u(3)],
  [r](:@{.}[u(3)])
  [r](:@{.}[u(3)])
  [r](:@{.}[u(3)])
  [r](:@{.}[u(3)])
  )
  "0"-[r]^{a}-[r]^{a}-[u]^{b}-[u]^{b}
  -[r]^{a}-[r]^{a}
  -[u]^{b}
  "0"-[u]^{b}-[r]^{a}-[r]^{a}-[r]^{a}
  -[u]^{b}
  -[u]^{b}
  -[r]^{a}
}
$$
In a distributive lattice (and more generally in a modular lattice)
it is always the case that $j \jmup m$ if and only if $m \mjdown j$.
Indeed, if $j \jmup m$, then $m^{\star} = j \vee m$.  Since $m \land j
< j$, $m \land j \leq j^{\star}$ and if this inequality is strict,
then $j,j_{\star},m,m^{\star},j\land m$ form a pentagon. Recalling
that a distributive lattice is a modular semidistributive lattice, we
obtain the following well known consequences:
\begin{Lemma}
  If $L$ is a modular lattice, then the reflexive and transitive
  closure $\dreftrclosure$ of the join dependency relation is an
  equivalence relation. Consequently, the congruence lattice of $L$ is
  a Boolean algebra.
\end{Lemma}
\begin{Lemma}
  \label{lemma:Ddistr}
  If $L$ is a distributive lattice, then the join dependency relation,
  is empty and there is a bijection between congruences on $L$ and
  subsets of join irreducible elements of $L$. 
\end{Lemma}

We illustrate congruences on the lattice path $\LP(n,m)$, since they
have a strong geometrical appealing.  Let us identify the join
irreducible element $\ji{x,y}$ with a disk within the interior of the
square $[x,x+1]\times [y-1,y]$: such a disk represents a hole in the
square $[0,n]\times [0,m]$ whose goal is to separate paths from
$(0,0)$ to $(n,m)$.  Let us fix a set $S$ of join irreducible elements
or holes. The formula \eqref{eq:congrfromS} that extracts the
congruence $\theta$ from the set $S$ can be interpreted by saying that
two paths $f,g$ are equivalent if and only if there is no hole in $S$
separating them. In the following diagram, we consider $\LP(3,3)$, we
let $S = \set{(0,3),(1,2)} $, and draw the resulting three equivalence
classes as paths up to dihomotopies \cite{goubault}.
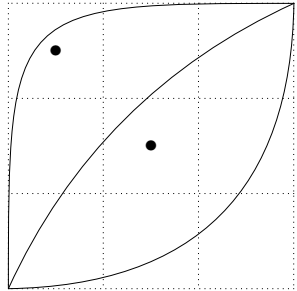
$$
\xygraph{
  []!M
  [u(3)][r(0.5)d(0.5)]*+{\bullet}
  "0"
  [r(1)u(2)][r(0.5)d(0.5)]*+{\bullet}
  "0":@{-}@`{"N","N"}"NE"
  "0":@{-}@/^5mm/"NE"
  "0":@{-}@/_13mm/"NE"
}
$$

\subsection{Lattices of permutations}
We review some facts on lattices of permutations $\Perm(n)$, $n\geq
0$, usually named permutoedra. Elements of $\Perm(n)$ are permutations
on the set $\setof{n}$ and the order -- known as the weak Bruhat order
-- is defined by means of its covering relation: $\sigma \covers
\sigma'$ iff $\sigma ' = \sigma \circ \sigma^{i}$ and $\ell(\sigma) <
\ell(\sigma')$, where $\ell(\sigma)$ is the minimum number $\ell$ such
that $\sigma$ can be written as the product of $\ell$ exchanges,
$\sigma = \sigma^{j_{1}}\circ \ldots \circ \sigma^{j_{\ell}}$.
 
Let us define an \emph{inversion}%
\footnote{%
  Usually an inversion is the inversion of a given permutation.
  Here we shall use this name coherently with the usage.
}% 
or \emph{disagreement} as an (unordered) pair $\{a,b\} \subseteq
\setof{n}$. If $a < b$, then we denote the inversion $\{a,b\}$ by
$a\under b$.  The set of all inversions on the set $\setof{n}$ will be
denoted $\inversions_{n}$ or simply $\inversions$ if $n$ is
understood. For a permutation $\sigma$, define
\begin{align*}
  D(\sigma) & = \set{a\under b \,|\,\sigma^{-1}(a) > \sigma^{-1}(b)}\,,
  &
  A(\sigma) & = \set{a\under b \,|\,\sigma^{-1}(a) < \sigma^{-1}(b)}\,.
\end{align*}
The first is the set of inversion or disagreements of $\sigma$, the
latter is the set of its agreements.  It is well known that
\begin{align*}
  \sigma \leq \sigma' & \text{ iff } D(\sigma) \subseteq D(\sigma')
  \text{ iff } A(\sigma') \subseteq A(\sigma)\,.
\end{align*}
Let us say that a subset of inversions $X\subseteq \inversions$ is
\emph{closed} if
\begin{align*}
  a\under b \in X \text{ and } b\under c \in X & \text{ implies } a\under c \in X\,.
\end{align*}
$X\subseteq \inversions$ is \emph{open} if and only if it is the
complement of a closed, that is if and only if
\begin{align*}
  a < b < c \text{ and } a\under c \in X & \text{ implies } a\under b
  \in X \text{ or } b\under c \in X \,.
\end{align*}
Finally, $X$ is \emph{clopen} if it is closed and open. The following
Proposition was used in \cite{yanagimoto} to prove that the
weak Bruhat order gives rise to a lattice.
\begin{Proposition}
  \label{prop:clopens}
  A subset $X \subseteq \inversions$ is clopen if and
  only if $X = D(\sigma)$ for some permutation $\sigma$.
\end{Proposition}
\begin{longversion}
\begin{proof}
  It is easily verified that the set of disagreements of a permutation
  is clopen. We shall verify the converse by induction on the
  cardinality of the clopen $X$.
  
  Observe first that the empty set is the set of inversions of the
  identity permutation.
  
  Therefore let $X$ be non empty clopen. Since $X$ is open, we can
  chose an inversion of the form $i\under i+1 \in X$ and let $Y = X
  \setminus \set{i\under i+1}$.  Let $\sigma_{i}$ be the exchange
  $(i,i+1)$.  We claim that $\sigma_{i}(a) < \sigma_{i}(b)$ if
  $a\under b \in Y$ and that the set $\sigma_{i} Y =
  \set{\sigma_{i}(a) \under \sigma_{i}(b) \,|\,a\under b \in Y }$ is
  clopen.
  
  The first statement follows from the fact that the restriction of
  $\sigma_{i}$ to both $\set{1,\ldots ,n}\setminus\set{i}$ and
  $\set{1,\ldots ,n}\setminus\set{i+1}$ is order preserving. Also, $Y$
  is closed, and therefore $\sigma_{i}(Y)$ is closed: if $\sigma_{i}
  a\under \sigma_{i} b, \sigma_{i} b'\under \sigma_{i} c \in
  \sigma_{i} Y$ and $\sigma_{i}(b) = \sigma_{i}(b')$, then $b = b'$
  and $a\under b,b\under c \in Y$, hence $a\under c \in Y$ and
  $\sigma_{i} a \under \sigma_{i} c \in \sigma_{i} Y$.

  We verify that $\sigma_{i} Y$ is open. Let $\sigma_{i} a < b <
  \sigma_{i} c$ and $\sigma_{i} a\under \sigma_{i} c \in \sigma_{i}
  Y$, we distinguish three cases:
  \begin{itemize}
  \item If $\{ \sigma_{i} a , b , \sigma_{i} c\} \cap \{i,i+1\}$ has
    at most one element then $\sigma_{i}$ is order preserving on $\{
    \sigma_{i} a , b , \sigma_{i} c\}$. It follows that $a <
    \sigma_{i} b < c$ so that $a\under \sigma_{i} b \in Y$ or
    $\sigma_{i} b \under c \in Y$. Consequently $\sigma_{i} a\under b
    \in \sigma_{i} Y$ or $b \under \sigma_{i} c \in \sigma_{i} Y$.
  \item Suppose $\sigma_{i} a = i$ and $b = i + 1$: then $\sigma_{i}
    b\under a = i\under i+1, a\under c \in X$, and therefore
    $\sigma_{i} b \under c \in X$. As a consequence $\sigma_{i}
    b\under c \in Y$ and $b\under \sigma_{i} c \in \sigma_{i} Y$.
  \item If $b = i$ and $\sigma_{i} c = i +1$, then one argues by
    duality as in the previous case.
  \end{itemize}
  
  Once we have established that $\sigma_{i} Y$ is clopen, we can
  assume by induction that $\sigma_{i} Y = D(\tau)$ for some
  permutation $\tau$, and that $X$ is the set of inversions for the
  permutation $\sigma_{i} \circ \tau$.  \
\end{proof}
\end{longversion}
For $X \subseteq \inversions$, we shall denote by $l(X)$ its
closure. %
\begin{longversion}
  Evidently the closure of $X$ is calculated by iterating from
  $X$ the following operation : 
  \begin{align*}
    T(X) & = X \cup \{ \,a\under c\, |\, a\under b,b\under c \in X
    \,\}\,.
  \end{align*}
\end{longversion} 
With  $r(X)$ we shall denote the interior of $X$, defined by $r(X) =
\neg l(\neg X)$.
%% \begin{align*}
%%   r(X) & = \neg l(\neg X)\,.
%% \end{align*}
\begin{Lemma}
  The closure $l(X)$ of an open $X \subseteq \inversions$ is open.
  The interior $r(X)$ of a closed $X \subseteq \inversions$ is closed.
\end{Lemma}
\begin{longversion}
\begin{proof}
  It is enough to show that $T(X)$ is open whenever $X$ is open.  Let
  $a\under c \in T(X)$ and $a < b < c$. If $a\under c \in X$ then the
  result is obvious, so let us suppose that $a\under d, d\under c \in
  X$ for some $d \in \setof{n}$. If $b = d$ then the result is
  obvious, otherwise either $b < d$ or $d < b$. Let us suppose $b <
  d$: either $a\under b \in X \subseteq T(X)$, or $b\under d \in X$ in
  which case $b\under d, d \under c \in X$ implies $b\under c \in
  T(X)$. An analogous argument works if $d < b$.
  
  The second property follows by duality.  \
\end{proof} 
\end{longversion}
The previous Lemma leads to the following representation of the
permutoedron, see  \cite{bjorner,barr}, and to simple formulas to
compute in the permutoedron. We shall often make use of this
representation  later.

Let $L$ be the Boolean algebra
of subsets of $\inversions$, $L_{r}$ be collection of open
subsets, and $L_{l}$ be collection of closed subsets and call $L_{rl}$
the collection of all clopens:
$$
\xygraph{
  []*+{L_{r}}="0"
  (:@{ >->}[r]*+{L}
  :@{ ->>}[d]*+{L_{l}}="END"^{l},
  :@{ ->>}[d]*+{L_{rl}}^{l}:@{ >->}"END"
  )
  [d(0.5)r(0.5)]*+{\vee}
  "0"[r(3)]*+{L_{l}}
  (:@{ >->}[r]*+{L}
  :@{ ->>}[d]*+{L_{r}}="END"^{r},
  :@{ ->>}[d]*+{L_{lr}}^{r}:@{ >->}"END"
  )
  [d(0.5)r(0.5)]*+{\land}
}
$$
The diagram on the left is meant to show that $L_{r}$ is a
sub-join-semilattice of $L$ and $L_{l}$ is a quotient-join-semilattice
of $L$; $L_{rl}$ a join-quotient of $L_{r}$, and a
sub-join-semilattice of $L_{l}$.  We obtain an useful formula for
computing the join of two clopens in $L_{rl}$: $X \vee Y = l(X \cup
Y)$.  The diagram on the right is meant to exemplify the dual notions,
with $L_{rl}$ replaced by $L_{lr}$ and join homomorphism replaced by
meet-homomorphism.  We obtain an useful formula for computing the meet
of two clopens in $L_{lr}$: $X \land Y = r(X \cap Y)$.  However
$L_{rl} = L_{lr}$ is simply
the collection of clopens, and in both cases the order is subset
inclusion: therefore we have $L_{lr} = L_{rl} $ as lattices.  

These formulas provide a method for computing meets and joins of
permutations, given in set-of-inversions form.  An explicit proof of
Proposition \ref{prop:clopens} suggests how to recover the string
representation of a permutation from its set-of-inversions
representation.  We recall that efficient algorithms for computing the
meet and the join of two permutations given in string representation
were proposed in \cite{markowsky}.

%%% Local Variables: 
%%% mode: latex
%%% TeX-master: "0"
%%% ispell-dictionary : "english"
%%% End: 

\section[The Lattices Structure of a Set of Multipermutations]%
{The Lattices Structure \\%
\mbox{\hspace{2cm}}of a Set of Multipermutations}
 
We shall consider paths in the space $\N^{n}$ from $0$ to a fixed
point $v$. These paths will have the property that each time step
increases just one coordinate. We  will denote these paths
by words over a totally ordered alphabet
$\Sigma = \set{a_{1},a_{2},\ldots, a_{n} }$%
\footnote{%
  We shall often make implicit the assumption that $\Sigma =
  \set{1,\ldots ,n}$.  } %
of directions.  If $v = (v_{1},\ldots ,v_{n})$ and $k = v_{1} + \ldots
+ v_{n }$, then we define
\begin{align*}
  \LP(v) & = \set{ w \in
    \Sigma^{k}\,|\,|w|_{a_{i}} = v_{i}, \text{ for } i =1,\ldots,
    n}\,.
\end{align*}
The set $\LP(v)$ is the set of multipermutations on $v$.  The
bijection between multipermutations in $\LP(v)$ and the paths we are
considering takes a $w \in \LP(v)$ to the path $f_{w}: \set{0,\ldots
  ,k} \rTo \N^{n}$ defined by $f_{w}(0) = 0$ and $f_{w}(i) =
f_{w}(i-1) + e_{w_{i}}$, where $e_{l} = (0,\ldots, 1,\ldots 0)$ has
just the coordinate $l$ different from $0$ and equal to $1$.

The rewrite  relation $\rightarrow$ on $\LP(v)$ is
defined  as follows:
\begin{align*}
  w \rightarrow u &
  \text{ iff } w = w_{1}a_{i}a_{j}w_{2},\;
  u = w_{1}a_{j}a_{i}w_{2},
   \text{ and } i < j\,.
\end{align*}
The rewrite relation is confluent and terminating, thus its reflexive
and transitive closure is a partial order $\leq$ on $\LP(v)$. W.r.t.
this partial order, the word $a_{1}^{v_{1}}a_{2}^{v_{2}}\ldots
a_{n}^{v_{n}}$ is the bottom element and
$a_{n}^{v_{n}}a_{n-1}^{v_{n-1}}\ldots a_{1}^{v_{1}}$ is the top.  We
have seen that in dimension 2 -- that is, for $n = 2$ -- this partial
order is a distributive lattice. In the general case $n \geq 3$, this
partial order is also a lattice, as we are going to argue.

It is harmless to assume that $v_{i} > 0$ for $i = 1,\ldots ,n$, so
that a word $w \in \LP(v)$ is a surjective
function $w : \setof{k} \rOnto \setof{n}$ %
                                %\footnote{Here we let $\setof{x} = \set{1,\ldots ,x}$} %
such that $|w^{-1}(i)| = v_{i}$ for $i = 1,\ldots ,n$.  Let $\mu^{\star} :
\setof{k} \rOnto \setof{n}$ be the only order preserving map with this
property and observe that for $w \in \LP(v)$ there exists a unique
permutation $\sigma : \setof{k} \rTo \setof{k}$ such that $w = \mu^{\star}
\circ \sigma$ and $\sigma$ is order preserving on every $w$-fiber.  We
shall denote this permutation by $\iota(w)$.  The following
Proposition was proved in \cite{bb} and  we are thankful to Peter
McNamara for independently pointing to us its proof.
\begin{Proposition}
  \label{prop:bijection}
  The function $\iota$ is an order isomorphism from the ordered set
  $(\LP(v), \leq)$ to a principal ideal of the permutoedron
  $\Perm(k)$.
\end{Proposition}
\begin{longversion}
\begin{proof}
  If $w = \mu^{\star} \circ \sigma$, then
  the condition that $\sigma$ is order preserving on every $w$-fiber
  -- that is, if $w_{i} = w_{j}$ and $i < j$, then $\sigma_{i} <
  \sigma_{j}$ -- can be rephrased in terms of agreements of $\sigma$:
  $k < l$ and $\mu^{\star}_{k} = \mu^{\star}_{l}$ implies $k\under l \in A(\sigma)$.
  Therefore the image of $\LP(v)$ under the correspondence $\iota$ is
  the set $I$ defined as
%% Since evidenly $\mu^{\star} \circ \sigma \in \LP(v)$, we have established a
%%   bijective correspondence between $\iota: \LP(v) \rTo I$, where
  \begin{align*}
    I 
    %% & = \set{ \sigma \in \Perm(k)\,| \,\forall k,l \;\; k<l \text{
    %%         and } \mu^{\star}_{k} = \mu^{\star}_{l} \text{ implies } k\under l \in
    %%       A(\sigma)
    %%     }\\
    & = \set{ \sigma \in \Perm(k)\,| \,\set{ k\under l \,|\, \mu^{\star}_{k} =
        \mu^{\star}_{l}} \subseteq A(\sigma) } \,.
  \end{align*}
  Since the set $\set{ k\under l \,|\, \mu^{\star}_{k} = \mu^{\star}_{l}}$ is clopen,
  $I$ is a principal ideal of $\Perm(k)$.
    
  Observe now that $w \rightarrow u$ if and only if $u = w \circ
  \sigma^{i}$ for some $i \in \set{1,\ldots ,k-1}$ such that $w_{i} <
  w_{i + 1}$. Hence, $\mu^{\star} \circ \iota(u) = \mu^{\star} \circ \iota(w) \circ
  \sigma^{i}$ and $\mu^{\star}(\iota(w)_{i}) < \mu^{\star}(\iota(w)_{i + 1})$. Since
  $\mu^{\star}$ is order preserving, we have $\iota(w)_{i} < \iota(w)_{i + 1}$
  and $\iota(w)_{i}\under \iota(w)_{i + 1} \in A(\iota(w))$.  As a
  consequence, $A(\iota(w) \circ \sigma^{i}) = A(\iota(w)) \setminus
  \set{\iota(w)_{i}\under \iota(w)_{i+1}}$.  Observing that
  $\mu^{\star}(\iota(w)_{i}) = w_{i} \neq w_{i + 1} = \mu^{\star}(\iota(w)_{i+1})$, we
  deduce that $\iota(w) \circ \sigma^{i}\in I$.  By uniqueness,
  $\iota(u) = \iota(w) \circ \sigma^{i}$.
  
  We have argued that $w \rightarrow u$ implies $\iota(w) \rightarrow
  \iota(u)$. The converse clearly holds, so that $\iota$ preserves and
  reflects the generating relation of the orders on $\LP(v)$ and
  $\Perm(k)$. It is therefore an order embedding.  
\end{proof}
\end{longversion}

Proposition \ref{prop:bijection} suggests that properties of $\LP(v)$
can be reduced to properties of the permutoedron.  For example:
\begin{Corollary}
  The lattice $\LP(v)$ is a bounded lattice.
\end{Corollary}
This follows from \cite{caspard} and \cite[\S 2.14]{freese}. Later,
our characterization of congruences of $\LP(v)$ will provide us with
another proof of this fact.  We can also use Proposition
\ref{prop:bijection} to argue that $w$ is join irreducible in $\LP(v)$
if and only if $\iota(w)$ is join irreducible in $\Perm(k)$.  It looks
unnatural, however, to deduce all the properties of $\LP(v)$ from the
representation above. For example, $\iota(w)$ need not to be meet
irreducible even if $w$ is such. Also, observe that permutoedra are
complemented lattices, while $\LP(v)$ is not: for $n = 2$, $\LP(v)$ is
a distributive lattice without necessarily being a Boolean algebra.
We can remark differences with lattices of paths in dimension 2 as
well, in particular $\LP(v)$ need not be distributive.  Therefore, we
seek for a direct understanding of $\LP(v)$, the key idea being the
equality
\begin{align*}
  \Perm(k) & = \LP(\underbrace{1,\ldots ,1}_{k-\text{times}})\,.
\end{align*}
We shall consider the lattice $\LP(v)$ as a generalization of the
permutoedron. The first step towards understanding its structure is to
find a working analogue of the notion of disagreement/agreement.

\begin{Definition}
  Given $w \in \LP(v)$ and $1 \leq l < m \leq n$, we define
  $\pi_{l,m}(w) \in \LP(v_{l},v_{m})$ as the word that arises by
  erasing all the symbols different from $a_{l}$ or $a_{m}$ (and by
  identifying the letter $a_{l}$ with $a$ and $a_{m}$ with $b$).
\end{Definition}
The formal definition of $\pi_{l,m}$ is that of a monoid morphism by
induction on the length of words.
\begin{Proposition}
  \label{lemma:projections}
  Let $w,u \in \LP(v)$, then $w \leq u$ if and only if $\pi_{l,m}(w)
  \leq \pi_{l,m}(u)$ for all $l,m$ such that $1 \leq l < m \leq n$.
\end{Proposition}
\begin{proof}
  Let us compute $\iota(w)$ for $w \in \LP(n,m)$: if $w_{\leq j}$ is
  the prefix of length $j$ of $w$, then
  \begin{align*}
    \iota(w)(j) & = 
    \begin{cases}
      |w_{\leq j}|_{a} & \text{ if } w_{j} = a \\
      n + |w_{\leq j}|_{b} & \text{ if } w_{j} = b\,.
    \end{cases}
  \end{align*}
  Clearly, $\iota(w)$ is a bijection, if we let $\mu^{\star}$ be the
  function sending $x$ to $a$ if $x \leq n$ and to $b$ otherwise then
  $\mu^{\star} \circ \iota(w) = w$, and finally if $i \leq j$ and
  $w_{i} = w_{j}$ then $\iota(w)(i) \leq \iota(w)(j)$.
  
  Observe that $\iota(w)^{-1}(i)$ is the length of least prefix of $w$
  containing $i$ $a$'s if $w_{i} = a$, and the length of least prefix of
  $w$ containing $i- n$ $b$'s if $w_{i} = b$.  We deduce that $i\under j
  \in D(\iota(w))$ iff the $(j - n)$-th occurrence of $b$ in $w$
  precedes the $i$-th occurrence of $a$ in $w$, and this happens if
  and only if the join irreducible $\ji{i-1,j}$ is below $w$. Taking
  into account this bijection between join irreducible elements below
  $w$ and inversions in $D(\iota(w))$, we conclude that $w \leq w'$ if
  and only if $D(\iota(w))\subseteq D(\iota(u))$.
  
  We consider now $w \in \LP(v)$ with $v \in \N^{n}$ and $n \geq 3$.
  Consider that the set $D(\iota(w))$ is the disjoint union of the
  sets
  \begin{align*}
    D_{l,m}(\iota(w)) & = \set{ i\under j \,|\, \sigma^{-1}(j) <
      \sigma^{-1}(j),\, \mu^{\star}(i) = l, \, \mu^{\star}(j) = m }\,,
  \end{align*}
  for $l,m$ such that $1 \leq l < m \leq n$.  For $i \in 1,\ldots ,n$
  let $k_{i} = \sum_{j = 1 \ldots l -1} v_{j}$.  Then $i\under j \in
  D_{l,m}(\iota(w))$ iff $i - k_{l}\under j - k_{m} + v_{l} \in
  D(\iota(\mu^{\star}_{l,m}(w))$ so that the two sets are in bijection.
  Therefore
  \begin{align*}
    w \leq u & \text{ iff } D(\iota(w)) \subseteq D(\iota(u))
    %\\& 
    \text{ iff } D_{l,m}(\iota(w)) \subseteq D_{l,m}(\iota(u)) 
    \tag*{whenever $1 \leq l <m \leq n$,}
  \end{align*}
  and, by the bijection, this holds iff $D(\iota(\pi_{l,m}(w)))
  \subseteq D_{l,m}(\iota(\pi_{l,m}(u)))$, that is $\pi_{l,m}(w) \leq
  \pi_{l,m}(u)$.  
\end{proof}

%%% Local Variables: 
%%% mode: latex
%%% TeX-master: "0"
%%% ispell-dictionary: "english"
%%% End: 

\section{Join and Meet Irreducible Elements in $\Paths(v)$}

A word $w \in \LP(v)$ is join irreducible iff there exists
a unique $i \in \set{1,\ldots ,n-1}$ such that $w_{i} > w_{i + 1}$.
Therefore we can write
\begin{align}
  \label{eq:joinirr}
  w & = (a_{1}^{x_{1}}a_{2}^{x_{2}}\ldots a_{n}^{x_{n}})
  (a_{1}^{\bar{x}_{1}}a_{2}^{\bar{x}_{2}}\ldots a_{n}^{\bar{x}_{n}})
\end{align}
where $\bar{x}_{i} = v_{i}- x_{i}$ for $i =1,\ldots ,n$. For such join
irreducible element $w$, we let $\vec{w}$ be the vector $x_{1},\ldots
,x_{n}$, so that $0 \leq \vec{w} \leq v$. Now let $x$ be any vector in
the closed interval $[0, v]$
and define
\begin{align*}
  \minji{x} & = \min \set{ i \,|\, x_{i} < v_{i} }&
  \maxji{x} & = \max \set{ i \,|\, x_{i} > 0 } 
\end{align*}
where for $x = v$, we let $\min \ji{x} = \infty$, and for $x = 0$ we
let $\max \ji{x} = - \infty$. If $x = \vec{w}$ for a join irreducible
element $w \in \Paths(v)$, then
\begin{align}
  \label{eq:condjoinirr}
  \minji{x} & <   \maxji{x}\,.
\end{align}
Every join irreducible is uniquely determined by a vector $0 \leq x
\leq v$ satisfying \eqref{eq:condjoinirr} and we shall use the
notation $\ji{x}$ for the word $w$ defined from $x$ in equation
\eqref{eq:joinirr}.
Observe that $\ji{x}$ can be defined for vectors $x$ for which
$\minji x \not< \maxji{x}$, in this case $\ji{x} = \bot$.  
\begin{Definition}
  We say that $(\minji{x}, \maxji{x})$ is the \emph{principal plan} of
  the join irreducible element $\ji{x}$.
\end{Definition}
Observe that a join irreducible element  $\ji{x}$ is uniquely
determined by the restriction of the vector $x$ to the closed interval
delimited by the principal plan: indeed,  $x_{i} = v_{i}$ for $i <
\minji{x}$, and  $x_{i} = 0$ for $i> \maxji{x}$.

\begin{shortversion}
  By counting vectors failing \eqref{eq:condjoinirr}, the following
  formula for the number of join irreducible elements in $\Paths(v)$
  was obtained in \cite{bb}:
  \begin{align*}
    \# \set{ w \in \Paths(v)\,|\,w \text{ is join irreducible}\,} & =
    \prod_{i =1 \ldots n} \!\!(v_{i} + 1) - (1 + \!\!\sum_{i = 1,\ldots
      ,n}\!\! v_{i})\,.
  \end{align*}
\end{shortversion}
\begin{longversion}
  Now consider a vector that fails \eqref{eq:condjoinirr}, i.e. for
  which $\maxji{x} \leq \minji{x}$: either $x_{i} = 0$ for $i =
  1,\ldots ,n$, or we can find an $i \in \set{1,\ldots ,n}$ such that
  $x_{i} > 0$ and $x_{j} = v_{j}$ for $j < i$ and $x_{j} = 0$ for $j >
  i$. By counting these vectors we obtain a formula for the number of
  join irreducible elements in $\Paths(v)$:
\begin{align*}
  \# \set{ w \in \Paths(v)\,|\,w \text{ is join irreducible}\,} & =
  \prod_{i =1 \ldots n} \!\!(v_{i} + 1) - (1 + \!\!\sum_{i = 1,\ldots
    ,n}\!\! v_{i})\,.
\end{align*}
\end{longversion} 
Analogous considerations hold for meet irreducible elements: a $w \in
\LP(v)$ is meet irreducible iff there exists a unique $i \in
\set{1,\ldots ,n-1}$ such that $w_{i} < w_{i + 1}$.  We can write
\begin{align}
  \label{eq:meetirr}
  w & = (a_{n}^{x_{n}}a_{n-1}^{x_{n-1}}\ldots a_{1}^{x_{1}})
  (a_{n}^{\bar{x}_{}}a_{n-1}^{\bar{x}_{n-1}}\ldots
  a_{1}^{\bar{x}_{1}})\,,
\end{align}
and defining
\begin{align*}
  \minmi{x} & = \min \set{ i \,|\, x_{i} > 0 }\,, &
  \maxmi{x} & = \max \set{ i \,|\, x_{i} < v_{i} }\,,
\end{align*}
where for $x = 0$ we let $\minmi{x} = \infty$ and for $x = v$ we let
$\maxmi{x} = -\infty$, we observe that $\minmi{x}<\maxmi{x}$.  Every
meet irreducible is uniquely determined by such a vector $0 \leq x
\leq v$ and we use the notation $\mi{x}$ for the $w$ of
\eqref{eq:meetirr}.  Again $\mi{x}$ is well defined even if $x$ does
not satisfy $\minmi{x} < \maxmi{x}$, in this case $\mi{x} = \top$.
$(\minmi{x}, \maxmi{x})$ is the principal plan of a meet
irreducible $\mi{x}$.

\begin{Lemma}
  \label{lemma:mainplan}
  Let  $\ji{x} \in \LP(v)$ be join irreducible, and let $(m,M)$
% = (\min  \ji{x},\max \ji{x})$
  be its principal plan. Then:
  \begin{itemize}
  \item  either $\pi_{i,j}(\ji{x})$ is
    join irreducible, or it is $\bot$.
  \item $\pi_{m,M}(\ji{x})$ is join irreducible,
  \item if $\pi_{i,j}(\ji{x})$ is  join irreducible, then $[i,j]
    \subseteq [m,M]$.
  \end{itemize}
\end{Lemma}
\begin{proof}
  If $w$ is the join irreducible of \eqref{eq:joinirr}, then
  $\pi_{i,j}(w) =
  a_{i}^{x_{i}}a_{j}^{x_{j}}a_{i}^{\bar{x}_{i}}a_{j}^{\bar{x}_{j}}$,
  and therefore $\pi_{i,j}(\ji{x}) = \ji{x_{i},x_{j}}$.

  We have $\pi_{m,M}(\ji{x}) = \ji{x_{m},x_{M}} $ and by definition 
  $x_{m} < v_{m}$ and $0 < x_{M}$.
  
  If $\ji{x_{i},x_{j}} = \pi_{i,j}(\ji{x})$ is join irreducible, then
  $x_{i} < v_{i}$ and $0 <x_{j}$, and hence $m \leq i$ and $j \leq M$.
\end{proof}

We shall use the characterization of the unique join irreducible
$\kappa^{d}(\mi{y})$ such that $\mi{y} \mjdown \kappa^{d}(\mi{y})\jmup
\mi{y}$ in the distributive lattice $\LP(v_{i},v_{j})$, cf.
\eqref{eq:jmupllp}, to characterize the relation $\ji{x} \jmup \mi{y}$
in $\LP(v_{1},\ldots ,v_{n})$, $n \geq 3$. We shall use the notation
$x_{|(a,b)}$ for the restriction of the function/vector $x$ to the
open interval $(a,b)$. Hence $x_{(a,b)} = x'_{(c,d)}$ iff $x_{i} =
x'_{i}$ for all $i$ such that $a <i < b$. 
\begin{Proposition}
  \label{prop:jmuplp}
  Let $\ji{x}$ be join irreducible, $\mi{y}$ be meet irreducible, and
  let $(a,b)$ and $(c,d)$ be their respective principal plans.  The
  relation $\ji{x} \jmup \mi{y}$ holds if and only if
  $\pi_{c,d}(\ji{x}) = \kappa^{d}(\pi_{c,d}(\mi{y}))$ and $x_{|(c,d)}
  = y_{|(c,d)}$.
\end{Proposition}
\begin{proof}
  As a first step we claim that $\ji{x} \jmup \mi{y}$ iff
  $\pi_{c,d}(\ji{x}) = \kappa^{d}(\pi_{c,d}(\mi{y}))$ and
  $\pi_{i,l}(\ji{x}) \leq \pi_{i,l}(\mi{y})$ for $(i,l) \neq (c,d)$.
  
  Observe that $\pi_{c,d}(\mi{y}^{\star}) = \pi_{c,d}(\mi{y})^{\star}$
  and $\pi_{i,l}(\mi{y}^{\star}) = \pi_{i,l}(\mi{y})$ if $(i,l) \neq
  (c,d)$.  Therefore, the relation $\ji{x} \leq \mi{y}^{\star}$
  implies $\pi_{c,d}(\ji{x}) \leq \pi_{c,d}(\mi{y})^{\star}$; if
  $\pi_{c,d}(\ji{x}) \leq \pi_{c,d}(\mi{y})$ we would have $\ji{x}
  \leq \mi{y}$; hence $\pi_{c,d}(\ji{x}) \not\leq
  \pi_{c,d}(\mi{y})$; overall we obtain $\pi_{c,d}(\ji{x}) =
  \kappa^{d}(\mi{y})$, i.e.  $x_{c} = y_{c} -1$ and $x_{d} = y_{d} +
  1$.
  
  As a second step, we claim that the condition $\pi_{c,d}(\ji{x}) =
  \kappa^{d}(\pi_{c,d}(\mi{y}))$ and $(i,l) \neq (c,d)$ implies
  $\pi_{i,l}(\ji{x}) \leq \pi_{i,l}(\mi{y})$ is equivalent to the
  condition $\pi_{c,d}(\ji{x}) = \kappa^{d}(\pi_{c,d}(\mi{y}))$ and
  $x_{|(c,d)} = y_{|(c,d)}$.
  
  The condition is necessary.  Let $i \in (c,d)$, then the relation
  $\ji{x_{c},x_{i}} = \pi_{c,i}(\ji{x})\leq \pi_{c,i}(\mi{y}) =
  \mi{y_{c},y_{i}}$ holds and is equivalent to $y_{c} \leq x_{c} =
  y_{c} - 1$ or $x_{i} \leq y_{i}$: we deduce $x_{i} \leq y_{i}$.
  Similarly, we deduce $y_{i} \leq x_{i}$ from $\ji{x_{i},x_{d}}\leq
  \mi{y_{i},y_{d}}$ and therefore $y_{i} = x_{i}$.
  
  The condition is sufficient. We only need to prove that $(i,l) \neq
  (c,d)$ implies $\pi_{i,l}(\ji{x}) \leq \pi_{i,l}(\mi{y})$. If $i <
  c$ or $d < l$, then $\pi_{i,j}(\mi{y}) = \top$, see Lemma
  \ref{lemma:mainplan}. We suppose therefore that $c \leq i < l \leq
  d$ with $(i,j) \neq (c,d)$, for example with $c < i$: then $x_{i} =
  y_{i}$ which is enough to ensure the relation $\pi_{i,l}(\ji{x}) =
  \ji{x_{i},x_{l}} \leq \mi{y_{i},y_{l}} = \pi_{i,l}(\mi{y})$.  
\end{proof}

The following consequence of Proposition \ref{prop:jmuplp} is worth
remarking:
\begin{Corollary}
  Let $\ji{x}$ be join irreducible, $\mi{y}$ be meet irreducible, and
  let $(a,b)$ and $(c,d)$ be their respective principal plans.  If
  $\ji{x} \jmup \mi{y}$, then $a \leq c < d \leq b$.
\end{Corollary}
The principal plan of $\mi{y}$ is contained in the principal plan of
$\ji{x}$: the relation $\pi_{c,d}(\ji{x}) =
\kappa^{d}(\pi_{c,d}(\mi{y}))$ implies that $\pi_{c,d}(\ji{x})$ is
join irreducible, and therefore $[c,d] \subseteq [a,b]$ by Lemma
\ref{lemma:mainplan}. 

We rephrase explicitly Proposition \ref{prop:jmuplp} as follows:
\begin{align*}
  \ji{x} \jmup \mi{y} & \text{ if and only if } y_{c} = x_{c} + 1,
  y_{d} = x_{d} - 1, \text{ and }x_{i} = y_{i} \text{ for } i \in
  (c,d)\,.
\end{align*}
Dually,  $\mi{y} \mjdown \ji{x}$ holds if and only if
$\pi_{a,b}(\mi{y}) = \kappa(\pi_{c,d}(\ji{x}))$ and $y_{|(a,b)} =
x_{|(a,b)}$, i.e. $x_{a} = y_{a} - 1$,
$x_{b} = y_{b} + 1$,  and  $x_{i} = y_{i}$ for $i \in
(a,b)$.

\begin{Corollary}
  The lattice $\LP(v)$ is semidistributive.
\end{Corollary}
\begin{proof}
  If $\ji{x} \jmup \mi{y} \mjdown \ji{x}$, then $\ji{x}$ and $\mi{y}$
  have the same principal plan $(a,b)$, $y_{|(a,b)} = x_{|(a,b)}$,
  $y_{a} = x_{a} +1$ and $y_{b} = x_{b} -1$,  $y_{i} = 0$ for $i <
  a$, and $y_{i} = v_{i}$ for $i >b$.  These conditions uniquely
  determine a vector $y$ for which $\minmi{y} = a$ (since $y_{a} =
  x_{a} +1 > 0$ ) and similarly $\maxmi{y} = b$. Hence $\mi{y} =
  \kappa(\ji{x})$ is the unique meet irreducible with the property
  that $\ji{x} \jmup \mi{y} \mjdown \ji{x}$.
  
  Similarly, a join irreducible $\ji{x} = \kappa^{d}(\mi{y})$ such
  that $\mi{y} \mjdown \ji{x} \jmup \ji{y}$ is uniquely determined.
  
  It is easily verified that $\mi{y} = \kappa(\kappa^{d}(\mi{y}))$ and
  $\ji{x} = \kappa^{d}(\kappa(\ji{x}))$.  Therefore $\LP(v)$ is a
  semidistributive lattice.
\end{proof}

%% We seek now a characterization of the $D$-relation between join
%% irreducible elememts of $\LP(v)$.  To this goal, let us define
%% \begin{align*}
%%   \delta(n,m) & = 
%%   \begin{cases}
%%     1 & n < m \\
%%     0 & n = m \\
%%     -1 & n > m\,,
%%   \end{cases}
%% \end{align*}
%% where $n, m \in \N$.
%% For two join irreducible $\ji{x}, \ji{z}$ we
%% define $\delta(\ji{x}, \ji{z}) = (\delta(a,e),\delta(b,f))$ where
%% $(a,b)$ and $(e,f)$ are the principal plans of $\ji{x}$ and $\ji{z}$,
%% respectively.  Observe that if $[e,f] \subseteq [a,b]$, then
%% $\delta(\ji{x}, \ji{z}) \in \set{(0,0), (1,0),(0,-1),(1,-1)}$.

\begin{Definition}
  If $\ji{x}, \ji{z}$ are two join irreducible elements, of respective
  principal plans $(a,b)$ and $(e,f)$, then we say that $\ji{x} \dep
  \ji{z}$ if and only if
  \begin{itemize}
  \item $[e,f] \subseteq [a,b]$ and $z_{|(e,f)} = x_{|(e,f)}$,
  \item $z_{e} = x_{e} - d_{e}$ where $d_{e} \in \set{0,1}$ and $d_{e}
    = 0$ if $e = a$,
  \item $z_{f} = x_{f} + d_{f}$ where $d_{f} \in \set{0,1}$ and $d_{f}
    = 0$ if $f = b$.
  \end{itemize}
\end{Definition}

\begin{Proposition}
  \label{prop:D}
  For two join irreducible elements $\ji{x}, \ji{z}$, there exists a
  meet irreducible $\mi{y}$ such that $\ji{x} \jmup \mi{y}$ and
  $\mi{y} \mjdown \ji{z}$ if and only if $\ji{x} \dep \ji{z}$.
\end{Proposition}
\begin{proof}
  The condition is necessary.  Let $(c,d)$ the principal plan of
  $\mi{y}$, then $[e,f] \subseteq [c,d] \subseteq [a,b]$ and similarly
  \begin{align*}
    z_{|(e,f)} & = y_{|(e,f)} = (y_{|(a,b)})_{|(e,f)} 
    = (x_{|(a,b)})_{|(e,f)} = x_{|(e,f)}\,.    
  \end{align*}
  Let us consider $x_{e},y_{e},z_{e}$ and suppose first that $a < e$:
  if $c < e$, then $z_{e} = y_{e} - 1$ and $y_{e} = x_{e}$ imply
  $z_{e} = x_{e} - 1$, and if $c = e$ then $z_{e} = y_{e} - 1$ and
  $y_{e} = x_{e} + 1$ imply $z_{e} = x_{e}$.
  Similarly, if $e = c = a$, then $z_{e} = y_{e} - 1$ and
  $y_{e} = x_{e} + 1$ imply $z_{e} = x_{e}$.
  
  The condition is sufficient. To this goal, we need to define a
  vector $y$ such that $\ji{x} \jmup \mi{y} \mjdown \ji{z}$. If
  $(c,d)$ is the principal plan of $\mi{y}$, then $y$ is determined by
  the condition $y_{|(c,d)} = x_{|(c,d)}$, $y_{c} = x_{c} + 1$ and
  $y_{d} = x_{d} -1$. Thus we only need to define the principal plan
  $(c,d)$ of $\mi{y}$, which we do according to four possible cases:
  \begin{enumerate}
  \item $z_{e} = x_{e}$ and $z_{f} = x_{f}$: we let $(c,d) = (e,f)$,
  \item $a < e$ and $z_{e} = x_{e} -1$ and $z_{f} = x_{f}$: we let
    $(c,d) = (a,f)$, 
  \item $z_{e} = x_{e}$, $f < b$ and $z_{f} = x_{f} + 1$: we let
    $(c,d) = (e,b)$,
  \item $a < e < f < b$, $z_{e} = x_{e} -1$ and $z_{f} = x_{f} +1$:
    we let $(c,d) = (a,b)$. 
  \end{enumerate}
\end{proof}

Thus we see that $\ji{x}D\ji{z}$ if and only if $x \neq y$ and
$\ji{x}\dep\ji{z}$. The relation $\dep$ is clearly antisymmetric, from
which we see that if $\ji{x}D\ji{z}$, then $[e,f] \subset [a,b]$,
where $(a,b)$ and $(e,f)$ are the principal plans of $\ji{x}$ and
$\ji{z}$, respectively.  Therefore the $D$-relation contains no cycle
and by Lemma \ref{lemma:bounded} we obtain:
\begin{Corollary}
  $\LP(v)$ is a bounded lattice.
\end{Corollary}
%% \begin{proof}
%%   Recall that $\ji{x}D\ji{z}$ if and only if $\ji{x}\neq \ji{z}$ and
%%   $\ji{x} \jmup \mi{y} \mjdown \ji{z}$ for some $\mi{y}$.  We can see
%%   then that $\ji{x}D\ji{z}$ implies that $[e,f] \subset [a,b]$, where
%%   $(a,b)$ and $(e,f)$ are the principal plans of $\ji{x}$ and
%%   $\ji{z}$, respectively. Therefore the $D$-relation  contains no
%%   cycle and $\LP(v)$ is bounded by Lemma \ref{lemma:bounded}.  
%% \end{proof}

The relation $\dep$ is easily seen to be transitive.  Indeed, let
$\ji{z}$, $\ji{y}$, and $\ji{x}$ be three join irreducible elements,
with respective principal plans $(e,f)$, $(c,d)$ and $(a,b)$, and
suppose that $\ji{x}\dep \ji{y} \dep \ji{z}$. Clearly $z_{|(e,f)} =
x_{|(e,f)}$. If $a = e$, then $z_{e} = x_{e}$.  Suppose that $a < e$:
if $c < e$ then $y_{e} = x_{e}$ and $z_{e} = y_{e} - d_{e} = x_{e} -
d_{e}$ since $e \in (c,d)$ implies $y_{e} = x_{e}$, and if $c = e$,
then $z_{e} = y_{e} = x_{e} + d_{e}$, with $d_{e} \in \set{0,1}$.
Analogous considerations show that $z_{f} = x_{f} + d_{f}$ with
$d_{f}\in\set{0,1}$, hence $\ji{x}\dep  \ji{z}$.
\begin{Corollary}
  For two join irreducible elements $\ji{x}$ and $\ji{z}$, the pair
  $(\ji{x},\ji{x}_{\ast})$ belongs to the congruence
  $\theta(\ji{z},\ji{z}_{\ast})$ if and only $\ji{x}\dep\ji{z}$.
\end{Corollary}
Indeed, from what we have seen, the relation $\dep$ and the reflexive
and transitive closure $\dreftrclosure$ of the join dependency
relation coincide, and it is a general fact for finite lattices that
$(\ji{x},\ji{x}_{\ast}) \in \theta(\ji{x},\ji{x}_{\ast})$ if and only
$\ji{x}\dreftrclosure\ji{z}$.

It should also be observed that the explicit description of the
relation $\dreftrclosure$ suffices to compute the dimension monoid of
a lattice $\LP(v)$. According to \cite{wehrung} this is the
commutative monoid generated by join irreducible elements $j$ and
subject to the relations $j + k = j$ whenever $k \dreftrclosure j$.

\vspp

When $v = 1^{n}$, that is when $\LP(v) = \Perm(n)$, a join irreducible
element $\sigma$ is uniquely described by its principal plan $(a,b)$
and by the subset $D_{a}$ of the open interval $(a,b)$ of
disagreements of $a$, $D_{a} = \set{i \in (a,b)\,|\,a\under i \in
  D(\sigma)}$. In vector notation, if $\sigma = \ji{x}$, then $D_{a} =
\set{i \in (a,b)\,|\,x_{i} = 1}$. Taking the triple $(a,b,D_{a})$ as a
representation of a join irreducible, we have
\begin{Corollary}
  Let $(a,b,D_{a})$ and $(c,d,D_{c})$ be two distinct join irreducible
  elements of $\Perm(n)$. Then $(a,b,D_{a}) D (c,d,D_{c})$ if and only
  if $(c,d) \subseteq (a,b)$ and $D_{c} = D_{a} \cap (c,d) $.
\end{Corollary}

%% \begin{proof}
%%   Indeed, $(\ji{x},\ji{x}_{\ast}) \in \theta(\ji{x},\ji{x}_{\ast})$ if
%%   and only $\ji{x}D^{\ast}\ji{z}$, where $\dreftrclosure$ is the transitive
%%   closure of the join dependency relation.
  
%%   However the relation of Proposition \ref{prop:D}, out which the join
%%   dependency relation is defined, is already transitive. Let $\ji{z}$,
%%   $\ji{y}$, and $\ji{x}$ be three join irreducible elements satisfying
%%   the conditions of Proposition \ref{prop:D}, with respective
%%   principal plans $(e,f)$, $(c,d)$ and $(a,b)$. Clearly $z_{|(e,f)} =
%%   x_{|(e,f)}$. If $a = e$, then $z_{e} = x_{e}$.  Suppose that $a <
%%   e$: if $c < e$ then $y_{e} = x_{e}$ and $z_{e} = y_{e} - d_{e} =
%%   x_{e} - d_{e}$ since $e \in (c,d)$ implies $y_{e} = x_{e}$, and if
%%   $c = e$, then $z_{e} = y_{e} = x_{e} + d_{e}$, with $d_{e} \in
%%   \set{0,1}$. Analogous considerations show that $z_{f} = x_{f} +
%%   d_{f}$ with $d_{f}\in\set{0,1}$.  
%% \end{proof}

\vspp

Finally, for computational purposes, we study covers of the reflexive
transitive closure $\dreftrclosure$ of the join dependency relation.
In a semidistributive lattice, every cover of the $D$-relation is
either of type $A$ or of type $B$.  We recall that $j_{1}A j_{2}$ if
and only if $j_{1} \jmup \kappa(j_{2})$ and $j_{1} \neq j_{2}$, and
that $j_{1}B j_{2}$ if and only if $\kappa(j_{1}) \mjdown j_{2}$ and
$j_{1} \neq j_{2}$.  We refer the reader to \cite[\S 2.58]{freese} for
a general background on these relations. In the following Lemma
$d(x,y) = 1$ if $x \neq y$ and $d(x,y) = 0$ of $x = y$.
\begin{Lemma}
  Let $\ji{x},\ji{z}$ be join irreducible elements of $\LP(v)$, with
  respective principal plans $(a,b)$ and $(e,f)$. Then: 
  \begin{itemize}
  \item $\ji{x} A \ji{z}$ iff $(e,f) \subset (a,b)$ and $x_{|[c,d]} =
    y_{|[c,d]}$,
  \item $\ji{x} B \ji{z}$ iff $(e,f) \subset (a,b)$, $x_{|(c,d)} =
    y_{|(c,d)}$, $z_{e} = x_{e} - d(a,e)$, and $z_{f} = x_{f} +
    d(b,f)$.
  \end{itemize}
\end{Lemma}

It easily seen that if $\ji{x}D\ji{z}$ with $a < e$ and $f < b$, then
$\ji{x}D\ji{y}D\ji{z}$, where $\ji{y}$ is characterized by having
principal plan $(e,b)$ or by having principal plan $(a,f)$. We say in
the first case that $\ji{x}D\ji{y}$ is a left move, and in the second
case that it is a right move.  Say that the \emph{width} of a join
irreducible element is the distance between the two coordinates
forming the principal plan. Our next goal is to show that left moves
can be factorized through a sequence of left moves $\ji{x}D\ji{z}$
that decrease the width of the respective principal plans by one.
Clearly, an analogous result holds for right moves.
\begin{Lemma}
  Let $\ji{x}$, $\ji{z}$ be join irreducible elements of principal
  plans $(a,b)$ and $(e,b)$.  If $\ji{x}D\ji{z}$ with $a+ 1 < e$, then
  there exists a join irreducible element $\ji{y}$, of principal plan
  $(a+1,b)$, such that $\ji{x}D\ji{y}D\ji{z}$.
\end{Lemma}
\begin{proof}
  If $x_{a + 1} = v_{a + 1}$, then we let $y_{a + 1} = x_{a + 1} - 1$,
  otherwise, we let $y_{a + 1} = x_{a + 1}$. 
\end{proof}
\begin{Corollary}
  The set of join irreducible elements of $\LP(v)$ ordered by the
  relation $\dreftrclosure$ is a graded poset.
\end{Corollary}

\begin{Corollary}
  Every $D$-path of join irreducible elements in $\LP(v)$, $v \in
  \N^{n}$, has length at most $n-2$. If $v_{i} > 0$ for $i =
  1,\ldots ,n$, then such length is realized.
\end{Corollary}
\begin{proof}
  If $\ji{x} D \ji{y}$, then the principal plan of $\ji{y}$ is
  strictly contained in the principal plan of $\ji{x}$.  Conversely,
  consider the word $a_{n}^{v_{n}}a_{1}^{v_{1}}a_{2}^{v_{2}}\ldots
  a_{n-1}^{v^{n-1}}$. Permuting $a_{n}^{v^{n}}$ with $a_{i}^{v_{i}}$,
  $i = 1 \ldots n -2$ gives a $D$-chain of $n-2$ elements.  
\end{proof}

We accomplish our analysis of covers by classifying them in 4
categories. A lower cover $\ji{x}D\ji{z}$ can be $LA$, a left move of
type $A$, $LB$, a left move of type $B$, $RA$, a right move of type
$A$ and $RB$, a right move of type $B$.  We include next some
automatically generated examples.  The first diagram illustrates the
join dependency relation of $\Perm(4)$, which should be compared with
Figure 10 in \cite{duquenne}.
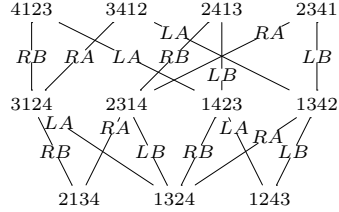
$$
\xygraph{%
%!{<1.5cm,0cm>:<0cm,1cm>::}%
[]="SLevel"%
[]*+{\scriptstyle{2134}}="2134" [r]*+{\scriptstyle{1324}}="1324" [r]*+{\scriptstyle{1243}}="1243" %
"SLevel"[u(1)l(0.5)]="SLevel"%
[]*+{\scriptstyle{3124}}="3124" [r]*+{\scriptstyle{2314}}="2314" [r]*+{\scriptstyle{1423}}="1423" [r]*+{\scriptstyle{1342}}="1342" %
"SLevel"[u]="SLevel"%
[]*+{\scriptstyle{4123}}="4123" [r]*+{\scriptstyle{3412}}="3412" [r]*+{\scriptstyle{2413}}="2413" [r]*+{\scriptstyle{2341}}="2341" %
"4123"(-"1423"|{LA},-"3124"|{RB})%
"3124"(-@/-2em/"1324"|{LA},-"2134"|{RB})%
"3412"(-@/-2em/"1342"|{LA},-"3124"|{RA})%
"2134"()%
"2413"(-@/0.7em/"1423"|{LB},-"2314"|{RB})%
"2314"(-"1324"|{LB},-@/-1em/"2134"|{RA})%
"2341"(-"1342"|{LB},-@/-2em/"2314"|{RA})%
"1423"(-@/-1em/"1243"|{LA},-"1324"|{RB})%
"1324"()%
"1342"(-"1243"|{LB},-@/-1em/"1324"|{RA})%
"1243"() %
}
$$

%%% Local Variables: 
%%% mode: latex
%%% TeX-master: "0"
%%% End: 

%
The following diagram illustrates the join dependency relations for
$\LP(2,1,1)$ and $\LP(1,2,1)$:
$$
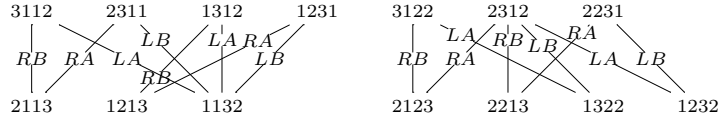
\xygraph{%
[]="SLevel"%
[]*+{\scriptstyle{2113}}="2113" [r]*+{\scriptstyle{1213}}="1213" [r]*+{\scriptstyle{1132}}="1132" %
"SLevel"[u]="SLevel"%
[]*+{\scriptstyle{3112}}="3112" [r]*+{\scriptstyle{2311}}="2311" [r]*+{\scriptstyle{1312}}="1312" [r]*+{\scriptstyle{1231}}="1231" %
"3112"(-"1132"|{LA},-"2113"|{RB})%
"2113"()%
"2311"(-@/-1em/"1132"|{LB},-"2113"|{RA})%
"1312"(-@/-0.7em/"1132"|{LA},-@/1em/"1213"|{RB})%
"1213"()%
"1231"(-"1132"|{LB},-@/-1.5em/"1213"|{RA})%
"1132"() %
}
\hspp
\xygraph{%
[]="SLevel"%
[]*+{\scriptstyle{2123}}="2123" [r]*+{\scriptstyle{2213}}="2213" [r]*+{\scriptstyle{1322}}="1322" [r]*+{\scriptstyle{1232}}="1232" %
"SLevel"[u]="SLevel"%
[]*+{\scriptstyle{3122}}="3122" [r]*+{\scriptstyle{2312}}="2312" [r]*+{\scriptstyle{2231}}="2231" %
"3122"(-@/-2em/"1322"|{LA},-"2123"|{RB})%
"2123"()%
"2312"(-"1232"|{LA},-@/-0.7em/"1322"|{LB},
-"2123"|{RA},-@/-0.7em/"2213"|{RB})%
"2213"()%
"2231"(-"1232"|{LB},-@/-1.4em/"2213"|{RA})%
"1322"()%
"1232"() %
}
$$

%%% Local Variables: 
%%% mode: latex
%%% TeX-master: "0"
%%% End: 

%%
The next diagram represents the $D$-relation in $\Paths(1,1,2,1,1)$:
$$
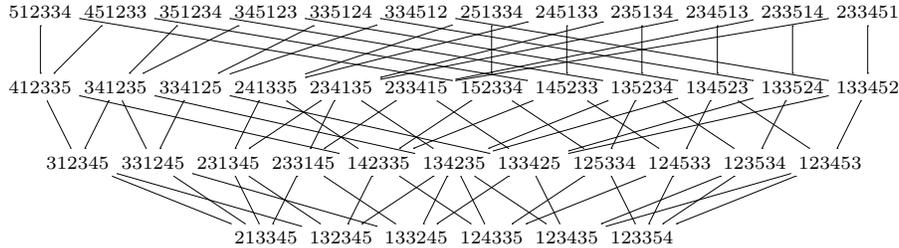
\xygraph{%
!{<1cm,0cm>:}
[]="SLevel"%
[]*+{\scriptstyle{213345}}="213345" [r]*+{\scriptstyle{132345}}="132345" [r]*+{\scriptstyle{133245}}="133245" [r]*+{\scriptstyle{124335}}="124335" [r]*+{\scriptstyle{123435}}="123435" [r]*+{\scriptstyle{123354}}="123354" %
"SLevel"[u(1)l(2.5)]="SLevel"%
[]*+{\scriptstyle{312345}}="312345" [r]*+{\scriptstyle{331245}}="331245" [r]*+{\scriptstyle{231345}}="231345" [r]*+{\scriptstyle{233145}}="233145" [r]*+{\scriptstyle{142335}}="142335" [r]*+{\scriptstyle{134235}}="134235" [r]*+{\scriptstyle{133425}}="133425" [r]*+{\scriptstyle{125334}}="125334" [r]*+{\scriptstyle{124533}}="124533" [r]*+{\scriptstyle{123534}}="123534" [r]*+{\scriptstyle{123453}}="123453" %
"SLevel"[u(1)l(0.5)]="SLevel"%
[]*+{\scriptstyle{412335}}="412335" [r]*+{\scriptstyle{341235}}="341235" [r]*+{\scriptstyle{334125}}="334125" [r]*+{\scriptstyle{241335}}="241335" [r]*+{\scriptstyle{234135}}="234135" [r]*+{\scriptstyle{233415}}="233415" [r]*+{\scriptstyle{152334}}="152334" [r]*+{\scriptstyle{145233}}="145233" [r]*+{\scriptstyle{135234}}="135234" [r]*+{\scriptstyle{134523}}="134523" [r]*+{\scriptstyle{133524}}="133524" [r]*+{\scriptstyle{133452}}="133452" %
"SLevel"[u]="SLevel"%
[]*+{\scriptstyle{512334}}="512334" [r]*+{\scriptstyle{451233}}="451233" [r]*+{\scriptstyle{351234}}="351234" [r]*+{\scriptstyle{345123}}="345123" [r]*+{\scriptstyle{335124}}="335124" [r]*+{\scriptstyle{334512}}="334512" [r]*+{\scriptstyle{251334}}="251334" [r]*+{\scriptstyle{245133}}="245133" [r]*+{\scriptstyle{235134}}="235134" [r]*+{\scriptstyle{234513}}="234513" [r]*+{\scriptstyle{233514}}="233514" [r]*+{\scriptstyle{233451}}="233451" %
"512334"(-"152334",-"412335")%
"412335"(-"142335",-"312345")%
"451233"(-"145233",-"412335")%
"312345"(-"132345",-"213345")%
"351234"(-"135234",-"341235")%
"341235"(-"134235",-"312345",-"331245")%
"345123"(-"134523",-"341235")%
"331245"(-"133245",-"213345")%
"335124"(-"133524",-"334125")%
"334125"(-"133425",-"331245")%
"334512"(-"133452",-"334125")%
"213345"()%
"251334"(-"152334",-"241335")%
"241335"(-"142335",-"231345")%
"245133"(-"145233",-"241335")%
"231345"(-"132345",-"213345")%
"235134"(-"135234",-"234135")%
"234135"(-"134235",-"231345",-"233145")%
"234513"(-"134523",-"234135")%
"233145"(-"133245",-"213345")%
"233514"(-"133524",-"233415")%
"233415"(-"133425",-"233145")%
"233451"(-"133452",-"233415")%
"152334"(-"125334",-"142335")%
"142335"(-"124335",-"132345")%
"145233"(-"124533",-"142335")%
"132345"()%
"135234"(-"123534",-"125334",-"134235")%
"134235"(-"123435",-"124335",-"132345",-"133245")%
"134523"(-"123453",-"124533",-"134235")%
"133245"()%
"133524"(-"123534",-"133425")%
"133425"(-"123435",-"133245")%
"133452"(-"123453",-"133425")%
"125334"(-"123354",-"124335")%
"124335"()%
"124533"(-"123354",-"124335")%
"123534"(-"123354",-"123435")%
"123435"()%
"123453"(-"123354",-"123435")%
"123354"() %
}
$$

%%% Local Variables: 
%%% mode: latex
%%% TeX-master: "0"
%%% End: 

%%
Finally, the join dependency relation in $\Paths(2,2,1,1)$ shows that
while a join irreducible element can have at most $4$ $D$-lower
covers, it may have more than 4 upper $D$-covers.
$$
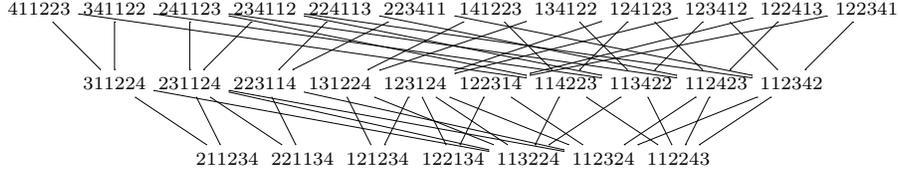
\xygraph{%
!{<1cm,0cm>:}
[]="SLevel"%
[]*+{\scriptstyle{211234}}="211234" [r]*+{\scriptstyle{221134}}="221134" [r]*+{\scriptstyle{121234}}="121234" [r]*+{\scriptstyle{122134}}="122134" [r]*+{\scriptstyle{113224}}="113224" [r]*+{\scriptstyle{112324}}="112324" [r]*+{\scriptstyle{112243}}="112243" %
"SLevel"[u(1)l(1.5)]="SLevel"%
[]*+{\scriptstyle{311224}}="311224" [r]*+{\scriptstyle{231124}}="231124" [r]*+{\scriptstyle{223114}}="223114" [r]*+{\scriptstyle{131224}}="131224" [r]*+{\scriptstyle{123124}}="123124" [r]*+{\scriptstyle{122314}}="122314" [r]*+{\scriptstyle{114223}}="114223" [r]*+{\scriptstyle{113422}}="113422" [r]*+{\scriptstyle{112423}}="112423" [r]*+{\scriptstyle{112342}}="112342" %
"SLevel"[u(1)l(1)]="SLevel"%
[]*+{\scriptstyle{411223}}="411223" [r]*+{\scriptstyle{341122}}="341122" [r]*+{\scriptstyle{241123}}="241123" [r]*+{\scriptstyle{234112}}="234112" [r]*+{\scriptstyle{224113}}="224113" [r]*+{\scriptstyle{223411}}="223411" [r]*+{\scriptstyle{141223}}="141223" [r]*+{\scriptstyle{134122}}="134122" [r]*+{\scriptstyle{124123}}="124123" [r]*+{\scriptstyle{123412}}="123412" [r]*+{\scriptstyle{122413}}="122413" [r]*+{\scriptstyle{122341}}="122341" %
"411223"(-"114223",-"311224")%
"311224"(-"113224",-"211234")%
"341122"(-"113422",-"311224")%
"211234"()%
"241123"(-"112423",-"114223",-"231124")%
"231124"(-"112324",-"113224",-"211234",-"221134")%
"234112"(-"112342",-"113422",-"231124")%
"221134"()%
"224113"(-"112423",-"223114")%
"223114"(-"112324",-"221134")%
"223411"(-"112342",-"223114")%
"141223"(-"114223",-"131224")%
"131224"(-"113224",-"121234")%
"134122"(-"113422",-"131224")%
"121234"()%
"124123"(-"112423",-"114223",-"123124")%
"123124"(-"112324",-"113224",-"121234",-"122134")%
"123412"(-"112342",-"113422",-"123124")%
"122134"()%
"122413"(-"112423",-"122314")%
"122314"(-"112324",-"122134")%
"122341"(-"112342",-"122314")%
"114223"(-"112243",-"113224")%
"113224"()%
"113422"(-"112243",-"113224")%
"112423"(-"112243",-"112324")%
"112324"()%
"112342"(-"112243",-"112324")%
"112243"() %
}
$$

%%% Local Variables: 
%%% mode: latex
%%% TeX-master: "0"
%%% End: 

%%
The exact number of upper $D$-covers clearly depends on the
multiplicities in the vector $v$.

%%% Local Variables: 
%%% mode: latex
%%% TeX-master: "0"
%%% ispell-dictionary: "english"
%%% End: 

\section{Dimension Equations for Multinomial Lattices}

\subsection{Pentagons}
For a quotient $a/b$ in a lattice we mean a pair of elements $a,b$
such that $a \leq b$. We shall say that a quotient is prime if $a
\prec b$, and write $a/b \subseteq c/d$ for $d \leq b \leq a \leq c$.

For a \emph{pentagon} in a lattice $L$ we mean a triple of elements
$a,b,c \in L$ such that
\begin{align*}
  b & \leq a \, &
  a \vee c & = a \vee b\,, &
  &\text{and} &
  a \land c & = a \land b\,.
\end{align*}
We denote such a triple by $N(a/b,c)$ and say that $a/b$ is the
central quotient of the pentagon. $N(a/b,c)$ is \emph{non degenerate}
if $b < a$.  In a pentagon $N(a/b,c)$, let $1_{N(a/b,c)} = a \vee c =
b \vee c$ and $0_{N(a/b,c)} = a \land c = b \land c$ (we shall omit
the subscripts if the underlying pentagon $N(a/b,c)$ is understood).
As usual, we say that two quotients $x/y$ and $z/w$ transposes to each
other if $x = y \vee z$ and $w = y \land z$ or viceversa. We denote
this relation by $\sim$, so that $x/y \sim z/w$ implies $\theta(x,y) =
\theta(z,w)$. In a pentagon $N(a/b,c)$ we have $c/0 \sim 1/b$ and $c/0
\sim 1/a$ (and dually $1/c \sim a/0$ and $1/c \sim 1/a$).  Hence, if a
pentagon $N(a/b,c)$ is non degenerate, then $\set{a,c}$ and
$\set{b,c}$ are antichains: if $c \leq a$ then $c = c \land a = 0$
hence $b \leq 1 = a$, and if $b \leq c$, then $c = b \vee c = 1$,
hence $b = 0 \leq a$.

\begin{Lemma}
  \label{prop:pentagon}
  \label{lemma:pentagon}
  Let $N(a/b,c)$ be a pentagon of a finite lattice and
  $a'/b' \subseteq a/b$ be a prime quotient.  Then we can find a prime
  quotient $x/y \subseteq b/0$ such that $(a,b) \in
  \theta(x,y)$.
\end{Lemma}
\begin{proof}
  We have $(a',b') \in \theta(a,0) = \theta(1,c) =
  \theta(b,0)$, hence $\theta(a',b') \subseteq
  \theta(b,0)$.  Moreover 
  \begin{align*}
    \theta(0,b) & = \bigvee_{b \geq x\succ
    y\geq 0} \theta(x,y)
  \end{align*}
  and since $\theta(a',b')$ is join prime in
  the lattice of congruences of $L$, $\theta(b',a') \subseteq
  \theta(x,y)$ -- that is, $a'/b' \in \theta(x,y)$ -- for a prime
  quotient $x/y \subseteq b/0$.  
\end{proof}
\begin{Lemma}
  Let $L$ be a lattice and $x/y$ be a prime quotient. Then we can find
  a join irreducible $j$ such that $j/j_{\star} \sim x/y$.
\end{Lemma}
\begin{proof}
  Consider the set of elements $z$ such that $z \vee y = x$. This set
  is non empty, let $j$ be minimal in this set. Observe
  that if $z$ is a lower cover of $j$, then $y \leq z \vee y \leq x$
  and, by minimality, $z \vee y = y$.
  It follows that $j$ is join irreducible, since if $z_{1},z_{2}$ are
  distinct lower covers, then $j \vee y = z_{1}\vee z_{2} \vee y = y$,
  contradicting $j \vee y =x$.
  
  Let $j_{\star}$ be the unique lower cover of $j$, we have seen that
  $j_{\star} \leq y$. It follows that
  $j_{\star} \leq j \land y \leq j$. Since
  $j \not\leq y$, it follows that $j_{\star} = j \land y$, showing that
  the covers $j/j_{\star}$ and $x/y$ transpose to each other. 
\end{proof}
%% Notice that if $L$ is join semidistributive -- i.e. it satisfies
%% \ref{eq:sdvee} only -- then the join irreducible $j$ in the
%% above statement is unique, hence we denote it by $j(x/y)$.

\begin{Corollary}
  \label{cor:pentD}
  If a lattice $L$ contains a non degenerate pentagon, then the $D$
  relation is non empty.
\end{Corollary}
\begin{proof}
  Let $N(a/b,c)$ be a non degenerate pentagon in $L$, chose a prime
  quotient $a'/b' \subseteq a/b$, and use Lemma \ref{lemma:pentagon}
  to find a prime quotient $x/y \subseteq b/0$ such that $(a',b') \in
  \theta(x,y)$. Hence, chose join irreducible elements $j,j'$ such
  that $j/j_{\ast} \sim a'/b'$ and $j'/j'_{\ast} \sim x/y$.  We have
  therefore $(j,j_{\ast}) \in \theta(j'/j'_{\ast})$, i.e.  $j
  \dreftrclosure j'$.  Notice finally that $j \neq j'$, since $j' \leq
  x \leq b'$ but $j \not\leq b'$, showing that $j \dtrclosure j'$,
  hence the join dependency relation $D$ is not empty.  
\end{proof}
Together with Lemma \ref{lemma:Ddistr}  we obtain a characterization
of distributive lattices.
\begin{Corollary}
  A finite lattice is distributive iff the join dependency relation,
  as well as the dual $D^{d}$ are empty.
\end{Corollary}
\begin{proof}
  Since distributivity is autodual, Lemma \ref{lemma:Ddistr} implies
  that $D$ and $D^{d}$ are empty in a distributive lattice.
  Conversely, observe that, given $j$, we can always find an $m$ and a
  $j'$ such that $j \jmup m \mjdown j'$. Then $j = j'$ since $D$ is
  empty, and such an $m$ is unique since $D^{d}$ is empty.  The given
  lattice is therefore semidistributive, modular by Corollary
  \ref{cor:pentD}, hence it is distributive.  
\end{proof}

\subsection{Meet semidistributivity at $n$ and the $D$-relation}
Recall that a (possibly infinite) lattice $L$ is meet semidistributive
if the relation  \ref{eq:sd}
%% \begin{align}
%%   \label{eq:sd}
%%   \tag{$SD(\land)$}
%%   x \land y  = x \land z & \Rightarrow x \land (y \vee z) = x\land y
%% \end{align}
holds in $L$. We are going to investigate such a condition within
finite lattices. To this goal, let $x,y,z$ be three variables, and
define two sequences of terms as follows:
\begin{align*}
  y_{0} & = y & z_{0} & = z \\
  y_{n + 1} & = y \vee (x \land z_{n} ) 
  & z_{n + 1} & = z \vee (x \land y_{n})\,.
\end{align*}
For each $n \geq 0$, the equation \ref{eq:sdn} is
\begin{align}
  \label{eq:sdn}
  \tag*{$SD_{n}(\land)$} x \land y_{n} & = x \land (y \vee z)\,.
\end{align}
If a lattice satisfies \ref{eq:sdn}, then we say that it is \emph{meet
  semidistributive at $n$}.  For an \emph{\ref{eq:sdn}-failure} we
mean a tuple $(L,x,y,z)$, where $L$ is a lattice, $x,y,z \in L$, and
$x \land y_{n} < x \land (y \vee z)$.\footnote{%
  Clearly, we are overloading notation, since we should make clear
  when we are dealing with terms, and when we are dealing with their
  interpretations in a given  lattice.  }

\vspp

We explicit some properties of the sequences $y_{n},z_{n}$.
\begin{Lemma} 
  For $n, m \geq 0$
  \begin{align*}
    y_{n} &\leq y_{n + 1}, & z_{n} & \leq z_{ n + 1},
    & \text{and } && y \vee z & =  y_{n} \vee z_{m}\,.
  \end{align*}
%%   $y_{n} \leq y_{n + 1}$ and $z_{n} \leq z_{ n + 1}$.
%%   For $m,n \geq 0$, $y \vee z =  y_{n} \vee z_{m}$, 
\end{Lemma}
\begin{proof}
  Clearly $y_{0} \leq y_{1}$ and $z_{0} \leq z_{1}$. If $y_{n} \leq
  y_{n + 1}$ and $z_{n} \leq z_{n + 1}$, then $y_{n+1} = y \vee
  (x\land z_{n}) \leq y \vee (x\land z_{n + 1}) = y_{n + 2}$, and
  similarly $z_{n + 1} \leq z_{n+ 2}$.
  
  For $y \vee z = y_{n} \vee z_{m}$, observe that the previous property
  implies $y \vee z \leq y_{n} \vee z_{m}$. Conversely, $y_{0} = y
  \leq y \vee z$, and $y_{n + 1} = y \vee (x\land z_{n}) \leq y \vee
  z_{n} \leq y \vee z$, assuming $z_{n} \leq y \vee z$.  Therefore
  $y_{n} \leq y \vee z$ and $z_{m} \leq y \vee z$ for each $n, m\geq
  0$, implying $y_{n} \vee z_{m} = y \vee z$.
\end{proof}
\begin{Fact}
  If \ref{eq:sdn} holds, then $SD_{k}(\land)$ holds for $k \geq n$.
\end{Fact}
%%\begin{proof}
Indeed, we can use the previous Lemma to see that if \ref{eq:sdn}
holds, then $ x \land (y \vee z) = x \land y_{n} \leq x \land y_{n +
  1} $ while $x \land y_{n + 1} \leq x \land ( y \vee z)$ holds for
each $n$.
%%\end{proof}

\begin{Fact}
  A finite lattice is meet semidistributive if and only if it is meet
  semidistributive at $n$ for some $n \geq 0$.
\end{Fact}
In \cite[\S 4.2]{jipsen} it is shown that a semidistributive
lattice generates a variety whose lattices are all meet
semidistributive if and only if one among the equations \ref{eq:sdn}
holds. It is straightforward to generalize the argument to obtain the
above statement. To this goal, 
let us define
\begin{align}
  \label{eq:xk}
  x_{k} & = (x \land y_{k-1}) \vee (x \land z_{k-1})\,,
  &&\text{for } k\geq 1. 
\end{align}

Let $x,y,z$ be fixed elements of a meet semidistributive lattice $L$,
we claim that if $y_{k} = y_{k + 1}$ and $z_{k} = z_{k + 1}$, then
$(L,x,y,z)$ is not an $SD_{k +1}(\land)$-failure.  Indeed
\begin{align*}
  y_{k} \land x & = y_{k+1} \land x \geq x_{k + 1} = (y_{k} \land x)
  \vee (y_{k} \land x) \geq y_{k} \land x \,.
\end{align*}
Therefore $y_{k} \land x = x_{k +1}$ and similarly
$z_{k} \land x  = x_{k +1}$. Hence
$y_{k} \land x = z_{k} \land x  =
x_{k +1}$ and therefore
\begin{align*}
 (y \vee z) \land x & = (y_{k} \vee z_{k}) \land x
 = x_{k+1} \leq y_{k+1}\,.
\end{align*}
Therefore if $L$ is finite and $x,y,z \in L$ we can let $\mu(x,y,z)$
be the least integer for which $y_{n-1} = y_{n}$ and $z_{n-1} =
z_{n}$. Then $L$ satisfies $SD_{M}(\land)$ where $M = \max
\set{\mu(x,y,z)\,|\, (x,y,z) \in L^{3}\,}$.

\vspp

The following Lemma shows that, when figuring out the configuration
given by $x$ and the sequences $y_{n}$ and $z_{n}$, we can assume that
$x \leq y \vee z$. While such representation has some interest for
heuristics, it plays no role in the following exposition.
\begin{Lemma}
  Let $x,y,z$ be given, let $\tilde{x} = x \land (y \vee z)$, and
  define the sequences $\tilde{y}_{n},\tilde{z}_{n}$ consequently out
  of the triple $\tilde{x},y,z$.  Then
  $\tilde{y}_{n} = y_{n}$, $\tilde{z}_{n} = z_{n}$, $x\land
  y_{n} = \tilde{x} \land y_{n}$, and $x\land z_{n} = \tilde{x} \land
  z_{n}$.
\end{Lemma}
\begin{proof}
  The relation $\tilde{y}_{n} = y_{n}$ and $\tilde{z}_{n} = z_{n}$
  hold for $n = 0$. By induction,
  \begin{align*}
    \tilde{x}_{n} \land \tilde{y}_{n} & = x \land (y \vee z) \land
    y_{n}
    = x \land
    y_{n}
  \end{align*}
  since $y_{n} \leq y \vee z$, hence $\tilde{z}_{n+ 1} = z_{n + 1}$.
\end{proof}

\begin{Proposition}
  Let $x_{k},k\geq1$ be defined as in $\eqref{eq:xk}$. Then
  \begin{align*}
 %%    x_{k} & = (x \land y_{k-1}) \vee (x \land z_{k-1})\,,
%%     &&\text{for } k\geq 1.
%%                                 %  \tag*{for $n\geq 1$.}
%%     \intertext{%
%%       Then%
%%     }
%%   %\begin{align*}
    x \land y_{k-1} & \leq x_{k} \leq x \land y_{k} \leq y_{k}
  \end{align*}
  and moreover $x \land y_{k}/x_{k}$ and $y_{k-1}$ form a possibly
  degenerate pentagon, with $x \land y_{k-1} = 0$ and $y_{k} = 1$.
\end{Proposition}
\begin{proof}
  Clearly $x \land y_{k-1} \leq x_{k}$ and $x_{k} \leq x$. To argue
  that $x_{k} \leq y_{k}$, observe that actually $y \vee x_{k} =
  y_{k}$:
  \begin{align*}
    y \vee x_{k} 
    & = y \vee (x \land y_{k-1}) \vee (x \land z_{k-1}) \\
    & = y  \vee (x \land z_{k-1}) \vee (x \land y_{k-1}) 
     = y_{k} \vee (x \land y_{k-1}) = y_{k}\,.
  \end{align*}
  It follows that $y_{k -1} \vee x_{k} =
  y_{k}$, and of course $(x \land y_{k}) \vee y_{k-1}  = x_{k} \vee
  y_{k-1} = y_{k}$, since $y_{k -1} \leq y_{k}$.
%%   \begin{align*}
%%     x \land y_{k} \land y_{k-1}
%%     & = x_{k} \land y_{k-1} =  x \land y_{k-1} \\
%%     (x \land y_{k}) \vee y_{k-1} & = x_{k} \vee y_{k-1} = y_{k}\,.
%%   \end{align*}
\end{proof}
Thus, for $k \geq 1$, we let $Y_{k} = N(x \land y_{k}/x_{k},y_{k
  -1})$. Since the roles of $y$ and $z$ are symmetric, $x \land
z_{k}/x_{k}$ and $z_{k-1}$ from a possibly degenerate pentagon which
we shall denote by $Z_{k}$.
 
The pentagons $Y_{k},Z_{k}$ might be degenerate, however, they respect
precise patterns: 
\begin{Lemma}
  \label{lemma:alternate}
  If $Z_{k}$ is non degenerate, then $Y_{k-1}$ is non degenerate.
\end{Lemma}
The Lemma is best proved with the following diagram at hand:
$$
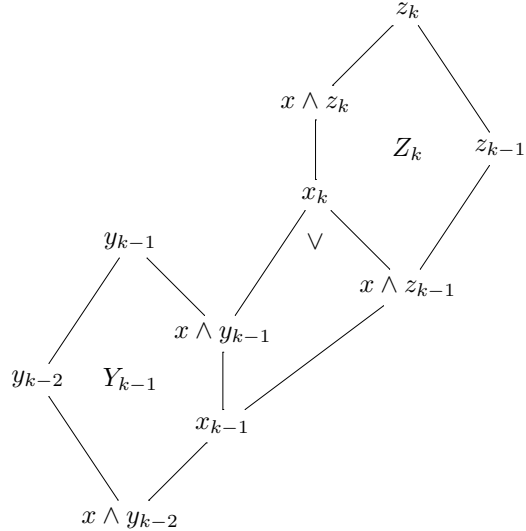
\mygraph[3.5em]{
  []!P{x_k}{x \land z_{k}}{z_{k -1}}{x \land z_{k -1}}{z_{k}}{Z_{k}}
  "1"="S2""0"="S1"
  "S"
  [l(3)d(2.5)]!Q{x_{k-1}}{x \land y_{k -1}}{y_{k -2}}{x \land y_{k
      -2}}{y_{k-1}}{Y_{k-1}}
  "2"="E2""1"="E1"
  "S1"-"E1"
  "S2"-"E2"
  "S2"[d(0.5)]*+{\vee}
}
$$
Suppose that $Y_{k-1}$ is degenerate, i.e. that $x_{k-1} = x \land
y_{k}$. Observe that 
$x_{k}/0_{Z_{k}}$ weakly projects down to $x \land y_{k-1}/x_{k-1}$,
i.e.
\begin{align*}
  x_{k} & = (x \land y_{k-1}) \vee (x \land z_{k-1})\,, \\
  x_{k-1} & = (x \land y_{k-2}) \vee (x \land z_{k-2}) \leq x \land
  z_{k-1} = 0_{Z_{k}}\,.
\end{align*}
Therefore, if the quotient $x \land y_{k}/x_{k-1}$ collapses, then
$x_{k} = 0_{Z_{k}}$. Since $x_{k}/0_{Z_{k}} \sim 1_{Z_{k}}/z_{k-1}$
and $1_{Z_{k}}/z_{k-1} \sim x \land z_{k}/0$, then $x \land z_{k} =
x_{k}$.
 
%% \begin{proof}
%% %  \myitem[\ref{item:successivi}]
%%   Let us suppose that $Y_{k}$ is degenerate, i.e.  $x_{k} = x \land
%%   y_{k}$.  We claim that $x_{k + 1} = x \land z_{k}$ and $z_{k + 1} =
%%   z_{k}$:
%%   \begin{align*}
%%     x_{k + 1} & = (x \land y_{k}) \vee (x \land z_{k})
%%     %\\ & 
%%     = x_{k} \vee (x \land z_{k}) 
%% %    & = (x \land y_{k-1}) \vee (x \land z_{k-1}) \vee (x \land z_{k}) \\
%%     = x \land z_{k}\,,
%%     \intertext{since $x_{k} \leq x \land z_{k}$, and } 
%%     z_{k + 1} & = z \vee (x \land y_{k})
%%     %\\ & = z \vee x_{k} 
%%     = z_{k}\,.
%%   \end{align*}
%%   Consequently $x_{k + 1} = x \land z_{k + 1}$ and $Z_{k + 1}$ is
%%   degenerate.  
%% \end{proof}

\begin{Lemma}
  Let $x,y,z$ be an \ref{eq:sdn}-failure in a meet semidistributive
  lattice $L$. Then,  either $Y_{k}$ is non degenerate, or
  $Z_{k}$ is non degenerate, for $1 \leq k \leq n$.
\end{Lemma}
\begin{proof}
%  \myitem[\ref{item:gemelli}]
  If $x_{k} = x \land y_{k}$ and $x_{k} = x \land z_{k}$, then $x
  \land y_{k} = x \land z_{k}$ implying $x \land (y \vee z) = x \land
  (z_{k} \vee y_{k}) = x \land y_{k}$. Therefore $SD_{k}(\land)$ holds and $k
  > n$.
  
  Consequently, if $1 \leq k \leq n$, then either $x_{k} < x \land
  y_{k}$, or $x_{k} < x \land z_{k}$.  
\end{proof}

\begin{Proposition}
  \label{prop:lenDpaths}
  Let $L$ be a finite meet semidistributive lattice such that
  \ref{eq:sdn} fails in $L$. Then $L$ contains a simple $D$-path of
  length $n$.
\end{Proposition}
\begin{proof}
  Let $(L,x,y,z)$ be an \ref{eq:sdn}-failure, and define the
  pentagons $Y_{k},Z_{k}$ as usual.  We have seen that either $Y_{n}$
  is non degenerate, or $Z_{n}$ is non degenerate. We shall suppose
  that $Y_{n}$ is non degenerate, so that $Y_{n - 2k}$ is non
  degenerate, and $Z_{n - 2k +1}$ is non degenerate, by Lemma
  \ref{lemma:alternate}.
  
  Since $Y_{n}$ is non degenerate, then we can chose a prime quotient
  $u_{n}/v_{n}$ contained in $x \land y_{n}/x_{n}$, the central
  quotient of the pentagon $Y_{n}$.

  Suppose that we have constructed a sequence of prime quotients 
  $$
  u_{n}/v_{n},\;  u_{n-1}/v_{n-1},\; \ldots ,\; u_{k}/v_{k}
  $$
  where $u_{i}/v_{i} \subseteq x \land y_{i}/x_{i}$ for $i$ even
  and otherwise $u_{i}/v_{i} \subseteq x \land z_{i}/x_{i}$; moreover
  $(u_{i}/v_{i}) \in \theta(u_{i - 1},v_{i - 1})$ for $i = n,\ldots
  ,k+1$.
  
  If $k > 0$, then we extend the sequence as follows. Since the roles
  of $Y_{i}$ and $Z_{j}$ are symmetric, we shall suppose that
  $u_{k}/v_{k}$ belongs to the central quotient of $Z_{k}$, i.e.
  $u_{k}/v_{k} \subseteq x\land z_{k}/x_{k}$. By Lemma
  \ref{lemma:pentagon} there is a prime quotient $u'/v' \subseteq
  x_{k} /x_{k} \land z_{k -1}$ such that $(u_{k},v_{k}) \in
  \theta(u'/v')$.
  
  If $k = 1$, then we let $u_{0}/v_{0} = u'/v'$.  Otherwise, recall
  that the quotient $x_{k}/x \land z_{k-1}$ weakly projects down to $x
  \land y_{k-1}/x_{k-1}$.  Consequently, $(u',v') \in \theta(x \land
  y_{k-1}/x_{k-1})$, hence -- as in Lemma \ref{lemma:pentagon} -- we
  can find a prime quotient $u_{k-1} /v_{k -1} \subseteq x \land
  y_{k-1}/x_{k-1}$ such that $(u',v') \in \theta (u_{k-1},v_{k-1})$.
  Since $(u_{k},v_{k}) \in \theta(u'/v')$, we have $(u_{k},v_{k}) \in
  \theta(u_{k -1}/v_{k -1})$ as well.
  
  Finally, for each prime quotient $(u_{i}/v_{i})$ let $j_{i}$ be a
  join irreducible element such that $j_{i} \vee v_{i} = x_{i}$ and
  ${j_{i}}_{\ast} \leq v_{i}$, so that $j_{i}/{j_{i}}_{\ast} \sim
  u_{i}/v_{i}$. Since $(j_{i},{j_{i}}_{\ast}) \in
  \theta(j_{i-1},{j_{i-11}}_{\ast})$, we have $j_{i} \dreftrclosure
  j_{i-1}$.
  
  Also observe that $j_{i}\neq j_{k}$ for $i \neq k$.  Indeed, let us
  suppose that $i < k$, then $j_{i} \leq u_{i} \leq (x\land y_{i})
  \vee (x \land y_{i}) = x_{i + 1} \leq x_{k}$ while $j_{k} \not\leq
  u_{k}$ implies $j_{k} \not\leq x_{k}$.  
  In particular we have $j_{i}  \dtrclosure
  j_{i-1}$ for $i = 1,\ldots ,n$. 
\end{proof}

\begin{Corollary}
  %\label{cor:Dinbounded}
  If $L$ is a meet distributive lattice whose maximal length of a
  simple $D$-path is $k$, then \ref{eq:sdn} holds in
  $L$ for $n > k$.
\end{Corollary}

\begin{Corollary}
  \label{cor:Dinbounded}
  If $L$ is a finite bounded lattice whose maximal length of a
  $D$-path is
  $k$, then \ref{eq:sdn} holds in $L$ for $n > k$.
\end{Corollary}
\begin{Corollary}
  For $v \in \N^{n}$, the multinomial lattice $\LP(v)$ is meet
  semidistributive at $k \geq n -1$.
\end{Corollary}
Indeed, the longest $D$-chain in $\LP(v)$ has length at most $n-2$. We
shall show later that this measure is tight, meaning that multinomial
lattices $\LP(v)$ fail to be meet semidistributive at $n -2$ provided
that $v_{i} > 0$ for $i = 1,\ldots ,n$.

\vspp

Finally, we exemplify the kind of problems arising when analysing
further Proposition \ref{prop:lenDpaths}.

It could seem natural generalize Corollary \ref{cor:Dinbounded} to non
bounded lattice by the following statement statement ``if every chain
of join irreducible elements in the congruence lattice of a meet
semidistributive $L$ has at most $k$ elements, then \ref{eq:sdn} holds
in $L$ for $n \geq k$''. The statement does not hold: let us consider
$\mathcal{JN}$, see \cite{jonssonnation}, the canonical example of a
non bounded semidistributive lattice, and its join dependency graph:
$$
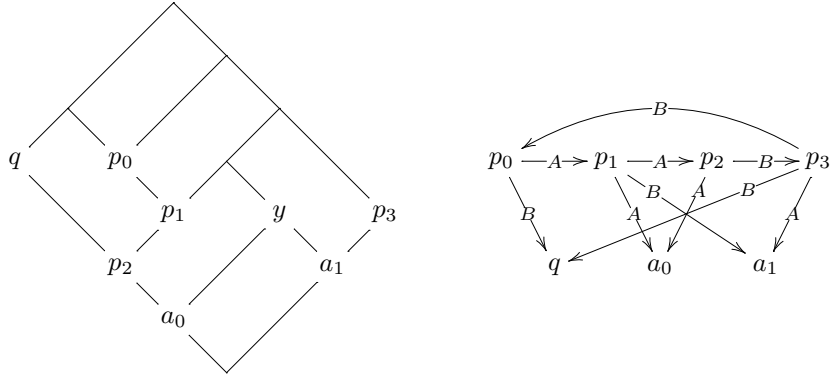
\mygraph[2.0em]{
  []="b"
  (-[lu]*+{a_{0}}="a0"
  -[lu]*+{p_{2}}="p2"(-[l(2)u(2)]*+{q}-[ru])
  -[ru]*+{p_{1}}="p1"(-[r(2)u(2)])
  -[lu]*+{p_{0}}="p0"(-[r(2)u(2)])
  -[lu]="c1"
  -[r(2)u(2)]="t"
  )
  (
  -[r(2)u(2)]*+{a_{1}}="a_{1}"
  (-[lu]*+{y}(-"a0")-[lu])
  -[ru]*+{p_{3}}="p_{3}"
  -[l(4)u(4)]
  )
}
\hspace{1cm}
\mygraph[4em]{%
  []*+{p_{0}}="p0"
  ([r]*+{p_{1}}="p1" [r]*+{p_{2}}="p2" [r]*+{p_{3}}="p3")
  [r(0.5)]([d]*+{q}="q"[r]*+{a_{0}}="a0"[r]*+{a_{1}}="a1")
  "p0"(:"p1"|{A},:"q"|{B})
  "p1"(:"p2"|{A},:"a0"|{A},:@/-1.5em/"a1"|{B})
  "p2"(:@/-1em/"a0"|{A},:"p3"|{B})
  "p3"(:"a1"|{A},:@/_2em/"p0"|{B},:@/-2.5em/"q"|{B})
}%
$$
%% In $\mathcal{JN}$ there are seven join irreducible elements, four
%% of which -- $p_{0},p_{1},p_{2},p_{3}$ -- lie on a cycle of the
%% $D$-relation . On the other hand, since an atom in a meet
%% semidistributive lattice is join prime, there is no $j$ such that
%% $a_{i}Dj$ for $i = 0,1$. Similarly, since a coatom in a join
%% semidistributive lattice is meet prime, if $c$ is a coatom then there
%% is no $m$ such that $cD^{d}m$. Since $\kappa(q)$ is a coatom, we
%% deduce that $q D j$ for no $j$.
It is easily argued that every
chain of join irreducible elements in the congruence lattice of
$\mathcal{JN}$ has at most 2 elements.  On the other hand,
$SD_{2}(\land)$ fails in $\mathcal{JN}$: let $x = p_{0}$, $z = q$, and
$y$ as in the diagram, then $x \land y_{2} = p_{1} < p_{0} = x = x
\land (y \vee z)$.

The following example exemplifies the problems found when trying to
find some converse to Proposition \ref{prop:lenDpaths}:
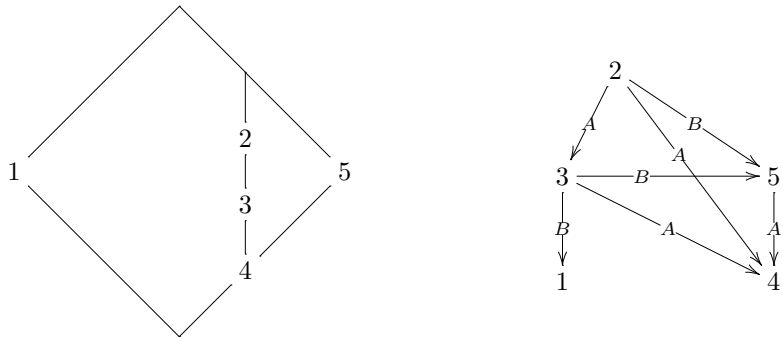
$$
\mygraph[2.5em]{
  []="b"
  (-[l(2.5)u(2.5)]*+{1}="1"
  -[r(2.5)u(2.5)]="t"
  )
  (
  -[ru]*+{4}="4"
  (-[u]*+{3}="3"-[u]*+{2}="2"-[u]="d"-"t")
  (-[r(1.5)u(1.5)]*+{5}="5"-"d")
  )
}
\hspace{2cm}
\mygraph[4em]{%
  []*+{2}="2"
  [l(0.5)]([d]*+{3}="3"[rr]*+{5})
  ([dd]*+{1}="1"[rr]*+{4})
  "2"(:"3"|{A},:@/-1em/"4"|{A},:"5"|{B})
  "3"(:"4"|{A},:@/-1em/"5"|{B},:"1"|{B})
  "5"(:"4"|{A})
}%
$$
Even if this bounded lattice contains simple $D$-paths of length
$3$ (and simple $A$-paths of length $2$) we cannot construct
\ref{eq:sdn} failures out of these paths. Actually, this lattice
 satisfies $SD_{2}(\land)$.

\subsection{Height of meet semidistributivity for $\Paths(v)$}

Let us say that the dimension of $\Paths(v)$ is the number  of indexes
$i$ such that $v_{i} > 0$. Every lattice $\Paths(v)$ of dimension $n$
is clearly isomorphic to a lattice $\LP(v')$ with $v' \in \N^{n}$.

\begin{Lemma}
  If $\LP(v)$ has dimension $n$, then $\Perm(n)$ embeds into $\LP(v)$.
\end{Lemma}
\begin{proof}
  We can assume that $v \in \N^{n}$, so that $v_{i}> 0$ for $i =
  0,\ldots ,n$.  Let $\psi : \Perm(n) \rTo \LP(v)$ such that
  $\psi(\sigma_{1}\ldots \sigma_{n}) = w_{\sigma_{1}}\ldots
  w_{\sigma_{n}}$ where $w_{j} = a^{v_{j}}_{j}$.
  
  Using Lemma \ref{lemma:projections} it is easily seen that this
  mapping is order preserving.  Now suppose that $\psi(\sigma) \vee
  \psi(\tau)\leq w$: therefore, if $i\under j$ is an inversion of
  either $\sigma$ or $\tau$, then all the $a_{j}$'s precede all the
  $a_{i}$'s in $w$.  By an easy induction, if $i\under j$ belongs to
  the closure of $I(\sigma) \cup I(\tau)$, then all the $a_{j}$'s
  precedes all the $a_{i}$'s in $w$, that is $\psi(\sigma \vee \tau)
  \leq w$.  Clearly, we can argue similarly to prove the relation
  $\psi(\psi \land \tau) = \psi(\sigma) \land \psi(\tau)$.  
\end{proof}

\begin{Proposition}
  If $n$ is the dimension of $\LP(v)$, then the equation $SD_{n
    -2}(\land)$ fails in $\LP(v)$.
\end{Proposition}
\begin{proof}
  Since $\Perm(n)$ is a sublattice of $\LP(v)$ it is enough to show
  that $SD_{n-2}(\land)$ fails in $\Perm(n)$. We shall use the
  representation of the permutoedron as a lattice of clopen sets.
  Hence, let us define
  \begin{align*}
    y &  = \set{ i \under i + 1 \,| \,i \text{ even}}\,, &
    z & = \set{i \under i + 1 \,| \,i \text{ odd} }\,, &
    x & = \set{1 \under i \, | \, i = 2,\ldots,n}\,.
  \end{align*}
  We claim that $x,y,z$ form an $SD_{n-2}(\land)$-failure in
  $\Perm(n)$. To ease the 
  verification, let
  $w_{0}  = \emptyset$ and 
  $w_{k} =  \set{1 \under i \, | \, i = 2,\ldots,k + 1}$,
  %%   \begin{align*}
  %%     w_{0} & = \emptyset\,, &
  %%     w_{k} & =  \set{1 \under i \, | \, i = 2,\ldots,k + 1},
  %%   \end{align*}
  where $1 \leq k \leq n - 1$, so that $x = w_{n - 1}$.  We claim
  that $w_{k} = x \land y_{k}$ if $k$ is even, and $w_{k} = x \land
  z_{k}$ if $w_{k}$ is odd.

  We remark first  that
  $x \land y_{0}  = \emptyset = w_{0}$ and 
  $x \land z_{1} = x \land z_{0}  = \set{1\under 2} = w_{1}$,
  %%   \begin{align*}
  %%     x \land y_{0} & = \emptyset = w_{0}
  %%     &
  %%     x \land z_{1} = x \land z_{0} & = \set{1\under 2} = w_{1}\,,
  %%   \end{align*}
  where we use the fact that $z_{0} \geq (x \land y_{0})$ implies that
  $z_{1} = z_{0}$.  We suppose therefore that $x \land y_{2k} =
  w_{2k}$ and $x \land z_{2k + 1} = w_{2k + 1}$.
  %%   \begin{align*}
  %%     x \land y_{2k} & = w_{2k} &
  %%     x \land z_{2k + 1} & = w_{2k + 1}\,.
  %%   \end{align*}
  We deduce
  \begin{align*}
    y_{2k + 2}  & = y \vee (x \land z_{2k + 1}) 
    = y \vee (w_{2k + 1}) 
    = l( y \cup w_{2k + 1})  \\
    & = l( \set{i \under i + 1 \,| \,i \text{ even} } \cup  \set{ 1
      \under i \, | \, i = 2,\ldots,2k + 2}) \\
    & = \set{i \under i + 1 \,| \,i \text{ even} } \cup  \set{ 1
      \under i \, | \, i = 2,\ldots,2k + 3}\,, 
     \\   
    \intertext{hence}
    x \land y_{2k + 2} & =
    r(\,\set{ 1
      \under i \, | \, i = 2,\ldots,2k + 3 }\,) \\
    & =
    \set{ 1
      \under i \, | \, i = 2,\ldots,2(k + 1) + 1} = w_{2(k + 1)}\,,
    \intertext{since this set is already open. Similarly:}
    z_{2(k + 1) + 1} & = z \vee (x \land y_{2k + 2}) 
    = z \vee w_{2k + 2} 
    = l( z \cup w_{2k + 2}) \\
    & = l( \set{i \under i + 1 \,| \,i \text{ odd} } \cup  \set{ 1
      \under i \, | \, i = 2,\ldots,2k + 3}) \\
    & = \set{i \under i + 1 \,| \,i \text{ odd} } \cup  \set{ 1
      \under i \, | \, i = 2,\ldots,2k + 4} \\
    \intertext{hence}
    x \land z_{2(k + 1) + 1} & =
    r(\,\set{ 1
      \under i \, | \, i = 2,\ldots,2(k + 1) +2 }\,) \\
    & =
    \set{ 1
      \under i \, | \, i = 2,\ldots,2(k + 1) + 2} = w_{2(k + 1) + 1}\,.
  \end{align*}

  Finally, observe that $y \vee z = \top$, hence $x \land (y \vee z) =
  x$.   Hence, if $n = 2k$ is even, then:
  \begin{align*}
    x \land y_{n-2} & = w_{n - 2} < w_{n - 1} = x \land (y
    \vee z)\,,
  \end{align*}
  and if $n = 2k + 1$ is odd, then  
  \begin{align*}
    x \land z_{n-2}  & = w_{n - 2} < w_{n - 1} = x \land (y
    \vee z)\,.
  \end{align*}
\end{proof}

%%% Local Variables: 
%%% mode: latex
%%% TeX-master: "0"
%%% ispell-dictionary : "english"
%%% End: 

%\bibliographystyle{latex8}
%\bibliography{biblio}

\end{document}